	\documentclass{article}
	\usepackage{amsthm,amsmath,amsfonts,amssymb,amsrefs,hyperref}
	\usepackage[T1]{fontenc}
	\usepackage{graphicx}
	\usepackage{color}
	\usepackage{hyperref}
	\usepackage{mathrsfs}
	\usepackage{dsfont}
	\usepackage[all]{xy} 

	\numberwithin{equation}{section}

	\newtheorem{theorem}{Theorem}[section]
	\newtheorem{definition}[theorem]{Definition}
	\newtheorem{proposition}[theorem]{Proposition}
	\newtheorem{lemma}[theorem]{Lemma}
	\newtheorem{corollary}[theorem]{Corollary}
	\newtheorem{remark}[theorem]{Remark}
	\newtheorem{example}[theorem]{Example}

	\def\C{\mathbb C}
	
	\def\K{\mathbb K}
	\def\N{\mathbb N}
	
	\def\R{\mathbb R}
	\def\Z{\mathbb Z}
	\def\Lsc{{\mathscr L}_{\textrm{sc}}(X)}
	\def\L2{L^{2}(\R^n)}
	\def\supp{\textrm{supp}~} 
	\def\XC{C^{\infty}_{\textrm{exp}}}

	\def\co{\colon}

	\def\ds{\displaystyle}
	\newcommand{\goesto}[2]{\overset{#1}{\underset{#2 \to \infty}{\longrightarrow}}} 

\begin{document}

	\title{\sc On the generation of groups of bounded linear operators on Fr\'{e}chet spaces}
	\author{
	  \small{\sc Arag\~ao-Costa, \'{E}der} R.\thanks{Partially supported by Grant: 2014/02899-3, S\~{a}o Paulo Research Foundation (FAPESP), Brazil.} \quad and \quad
	  \small{\sc Silva, Alex P.}\thanks{Supported by Grant: CAPES, Brazil.}\\
	\small University of S{\~a}o Paulo, \\
	\small 13566-590 S{\~a}o Carlos, SP, Brazil.
}

\maketitle

		
%


	\begin{abstract}

		In this paper we present a general method for generation of uniformly continuous groups on abstract Fr\'{e}chet spaces (without appealing to spectral theory) and apply it to a such space of distributions, namely ${\mathscr F}L^{2}_{loc}(\R^{n})$, so that the linear evolution problem
		\begin{equation*}
			\left\{\begin{array}{l}
			u_{t} = a(D)u, t \in \R \\
			u(0) = u_0
			\end{array}
			\right.
		\end{equation*} always has a unique solution in such a space, for every pseudodifferential operator $a(D)$ with constant coefficients. We also provide necessary and sufficient conditions so that the spaces $\L2$ and ${\mathscr E}'(\R)$ are left invariant by this group; and we conclude that the solution of the heat equation on ${\mathscr F}L^{2}_{loc}(\R^{n})$ for all $t \in \R$ extends the standard solution on Hilbert spaces for $t \geqslant 0$.

	
	\end{abstract}

\maketitle


\section{\bf{Introduction}}\label{Section:introduction}

\medskip

We consider problems of the type
	\begin{equation}\label{eq:Cauchy}
		\left\{\begin{array}{l}
			u' = Au, t \in \R \\
			u(0)=u_0 \in X
		\end{array}
		\right. ,
	\end{equation}
where $A \colon X \to X$ is a bounded linear operator on a Fr\'{e}chet space $X$. In order to solve it, we need to recognize the essence of the resolution on Banach spaces and adapt it to this more general formulation. In short, we will extend the usual results of generation of uniformly continuous groups of bounded linear operators on Banach spaces. See Pazy~\cite{Pazy}.

The main idea of this paper is based on a simple one, started out by Euler, Napier and Bernoulli's researches: the definition of the real exponential function $t\mapsto \exp(t)$,~\cites{Euler,Napier,Bernoulli}.

If there exists a differential function $u \colon (a,b) \to \R$ with the property that
	\[
		u'(t) = u(t), \mbox{ for every } t \in (a,b),
	\]
then we list some consequences:

	\begin{itemize}
	
	\item[i)] if $u(t_0)=0$ for some $t_0$ then $u \equiv 0$ in $(a,b)$;

	\item[ii)] there exists at most one function $u \colon \R \to \R$ which satisfies the differential problem
		\begin{equation}\label{eq:ODE-R}
			\left\{\begin{array}{l}
				u'(t) = u(t), t \in \R \\
				u(0) = 1
			\end{array}
			\right.
		\end{equation}
	
	\item[iii)] if $u$ satisfies \eqref{eq:ODE-R} then also satisfies
		\[
			u(t+s) = u(t)u(s) \mbox{ for every } s,t \in \R.
		\]

	\end{itemize}

	For every $t \in \R$, set
	\[
		u(t):=1+t+\dfrac{t^2}{2!}+\dfrac{t^3}{3!}+\cdots+\dfrac{t^n}{n!}+\cdots,
	\]
	which satisfies the Cauchy problem \eqref{eq:ODE-R}, since the sequence of the functions $(u_n)_{n \in \N}$ given by $u_1(t) := 1+t$ and $u_n(t) := u_{n-1}(t) + \dfrac{t^n}{n!}$ satisfies $u_n'=u_{n-1}$ and converges uniformly on intervals $[-M,M]$.

	Denote $u(t)$ by $\exp(t)$ or $e^{t}$, and set $e := u(1)$. Hence by definition, the map $t \mapsto e^{t}$ is the unique solution function of \eqref{eq:ODE-R}. More generally, for $a \neq 0$, the map $t \mapsto u(t) := u_0 e^{at}$ is the unique solution function of
		\begin{equation*}
			\left\{\begin{array}{l}
				u' = a u, t \in \R \\
				u(0)=u_0
			\end{array}
			\right. ,
		\end{equation*}
which is a particular case of \eqref{eq:Cauchy} with $X=\R$ and $A(y)=ay$ for $y \in \R$.

	Such construction may be extended in order to solve the Cauchy problems in $\R^N$ and more generally in infinite-dimensional Banach spaces. Indeed, given a bounded linear operator $A \colon X \to X$ on a Banach space $X=(X, \| \cdot \|_X)$, we define for every $t \in \R$ the bounded linear operator $\exp(tA) \colon X \to X$ by
		\[
			\exp(tA) := \sum_{n=0}^{\infty} \dfrac{t^n}{n!} A^n,
		\]
whence $t \mapsto u(t):= \exp(tA)u_0$ is the unique solution of \eqref{eq:Cauchy}.
	
The approach above is the core of linear (semi)groups definition and theory, and it has become a very useful tool for solving partial differential equations (PDEs) which involves a temporal derivative. Although it has been extensively studied over the last decades (the texts~\cites{Carvalho_book,Hale,Henry,Pazy,Yosida} are some of the most common), it is usually done on Banach spaces, which are not the most appropriate ones to deal with the distributional aspects involved.

We propose to extend such construction to more general spaces in order to solve a larger class of evolution problems. Shortly, the key point is the simple connection between the generator $A$ and the structure of the vector space $X$, so that the exponential of $A$ makes sense as a bounded linear operator $e^A \colon X \to X$ and may be used to solve the associated Cauchy problem.

This issue has already been dealt for other authors adding hypothesis on the generator or on the phase space $X$. In~\cites{Choe, Lemle_Wu,Kraaij}, the semigroup is assumed to be equicontinuous (in the sense of Banach-Steinhauss theorem for locally convex spaces, see~\cites{Narici,Osborne,Rudin}), while others researchers treat the question in some particular Fr\'{e}chet spaces, such as done by Dembart~\cite{Dembart} (who consider the phase space as the space of the continuous functions defined on $[a,b]$ into a fixed topological vector space $E$) and in~\cite{Leonhard} (setting $X=\K^{\N}$, that is, the collection of scalar sequences). In other words, none of them solves completely the problem of the (semi)group generation in the general case, so we feel confident to point out some results about it.

In order to recognize some of the advantages of this approach, we consider a typical linear problem over a Banach space and analyse it over another natural vector space, which is not normable, although it is a Fr\'{e}chet space.

\begin{example}

Let $C_B = C_B(\R, \C)$ be the vector space of the bounded uniformly continuous complex-valued functions on $\R$ with the supremum norm. For every $t \in \R$, set
	\[
		\big(T(t) \phi \big)(x):=\phi(t+x), \mbox{ for every } \phi \in C_B,
	\]
so that $\{T(t): t \in \R\}$ is a $C_0$-group of contractions on $C_B$ and its infinitesimal generator $A \colon D(A) \subset
C_B \to C_B$ (which exists since $C_B$ is a Banach space) is defined by
	\[
		\big(A \phi \big)(x) := \phi'(x)
	\]
for every $\phi \in D(A):=\{ \phi \in C_B: \phi' \mbox{ exists and belongs to } C_B\}$.

The closed densely defined operator $A$ is associated with the Cauchy problem
	\begin{equation*}
		\left\{\begin{array}{l}
			u_t = u_x, t \in \R \\
			u(0)=u_0 \in C_B
		\end{array}
		\right. 
	\end{equation*}
and its solution, namely $t \mapsto T(t)u_0$, is an infinitely differentiable function on $\R$ whenever $u_0 \in C_B \cap C^{\infty}(\R, \C)$; and these facts agree with Hille-Yosida theorem.

However we may consider this very same problem on $C^{\infty}=C^{\infty}(\R, \C)$, for which there is no norm which turns it into a Banach space. There is no difficulty on considering the family of operators $T(t)$ from $C^{\infty}$ into itself. Although $(t,x) \mapsto u_0(t+x)$ is a solution, we cannot explicit its generator by the usual theory, since the phase space is not a Banach space. See~\cite{Pazy,Yosida}.

Fortunately, $C^{\infty}$ is a Fr\'{e}chet space with the countable separating family of seminorms
	\[
		p_{(m,j)}(\phi) := \sup_{|x| \leqslant j} \left| \dfrac{d^m \phi}{dx^m}(x) \right|, \phi \in C^{\infty},
	\]
for $m\in \Z_+$ and $j \in \N$.

In section \ref{Section:applications}, we construct a dense subspace of $C^{\infty}$, namely $\XC = C^{\infty}_{\textrm{exp}}(\R, \C)$, where the exponential of the derivative operator is well defined. More precisely, we have that the partial sums $S_N := \ds\sum_{n=0}^{N} \dfrac{t^n}{n!}\dfrac{d^n}{dx^n}\phi$ converges in $C^{\infty}$ to a function in $\XC$, for every $\phi \in \XC$ and $t \in \R$; that is, if we set
	\[
		e^{t \frac{d}{dx}}\phi := \sum_{n=0}^{\infty} \dfrac{t^n \phi^{(n)}}{n!}, \mbox{ for every } \phi \in \XC,
	\]
then the operator $e^{t \frac{d}{dx}} \colon \XC \subset C^{\infty} \to \XC \subset C^{\infty}$ is well defined and is a bounded linear operator; and the family of operators $\{e^{t \frac{d}{dx}}: t \in \R\}$ is a uniformly continuous group on $\XC$ such that $\left( e^{t \frac{d}{dx}}\phi \right)(s) = \phi(s+t)$, for every $s \in \R$.

All these properties will be verified afterwards in Section \ref{Section:applications}. For now let us see some consequences of these properties. First, the derivative operator $A:=\frac{d}{dx} \colon \XC \to \XC$ is a bounded linear operator instead of a closed densely defined operator, whence it may be seen as the generator of the solution group. This information about the solution cannot be provided by the usual theory of Banach spaces.

If a function $\phi$ belongs to $\XC$ then it is an analytic function on $\R$. So in order to the series $\sum \frac{t^n \phi^{(n)}}{n!}$ be convergent in $C^{\infty}(\R)$, we restrict ourselves to functions $\phi$ extremely smooth; we mean, the set of the functions $\phi$ for which this series is convergent is very restrictive.

On the one hand, if we want to solve $(\ref{eq:Cauchy})$ for a linear operator more general than $d/dx$ and initial data $u_{0}$ not so smooth, the idea presented above does not fit anymore. However, we will present very simple conditions over a linear bounded operator on a Fr\'{e}chet space to ensure that its exponential exists and solves the problem \eqref{eq:Cauchy}, for every initial data on the space. Under such conditions, the notion of solution to the problem
		\begin{equation*}
			\left\{\begin{array}{l}
			u_{t} = a(D)u, t \in \R \\
			u(0) = u_0
			\end{array}
			\right.
		\end{equation*} 
will make sense, for every pseudodifferential operator $a(D)\colon {\mathscr S}(\R^{n}) \to {\mathscr S}(\R^{n})$ with constant coefficients.

\end{example}
		
\bigskip
This paper is organized as follows. In Section \ref{Section:preliminaries}, we provide a brief introduction to Fr\'{e}chet spaces and bounded linear operators on it; we present the main Fr\'{e}chet spaces we deal with, highlighting ${\mathscr F}L^{2}_{loc}(\R^{n})$ in Example \ref{FL2loc}; and define pseudodifferential operators on the Schwartz space, Definition \ref{Def.:pseudodiff}. The main references are~\cites{Folland,FollandPDE,Folland1992,HormanderI,HormanderIII,Narici,Osborne,Rudin,Stein,Treves,Treves1976,Yosida}. 

Then we discuss in Section \ref{Section:generation} the compatibility that a bounded linear operator is required to have in order to stablish the generation theorem (Theorem \ref{Th:generation}) on abstract Fr\'{e}chet space and to extend the theory of generation of uniformly continuous groups found in Pazy~\cite{Pazy} to this more general class of complete locally convex spaces. The requirements for the operator, found in Definition \ref{Def:LscX}, are simple ones and allows to achieve the results with small adaptations on the standard proofs.

Section \ref{Section:applications} is completely concerned with applying the results in Section \ref{Section:generation} to PDEs. By Theorem \ref{Th:generationFL2loc}, pseudodifferential operators with constant coefficients defined on ${\mathscr F}L^{2}_{loc}(\R^{n})$ satisfy the conditions of the generation theorem and then solve Cauchy problems for which the initial data belongs to this Fr\'{e}chet space of distributions. Also, Theorems \ref{Th:restriction-E-R} and \ref{Th:restriction-L2} stablish necessary and sufficient conditions under which the subspaces ${\mathscr E}'$ and $L^2$ of ${\mathscr F}L^{2}_{loc}$ are left invariant by such group, respectively. The Example \ref{Ex.:heat} analyse the heat equation and the relation between the standard solution on Hilbert spaces for $t\geqslant0$ and the distributional solution on ${\mathscr F}L^{2}_{loc}(\R^{n})$ for all $t \in \R$.

At last in Section \ref{Section:final}, we wonder about the several implications of this theory on our understanding of evolution problems and propose to study whether the spectral theory may be adapted in order to extend the Hille-Yosida theorem and others to Fr\'{e}chet spaces.


\section{\bf{Preliminaries}}\label{Section:preliminaries}

Before introducing our results, for the sake of understanding, we are going to recall some definitions and to point out some differences between Banach spaces and Fr\'{e}chet spaces.

	\begin{definition}
		A topological vector space (or a TVS, for short) is a vector space $X$ over a field $\K$ ($\R$ or $\C$) endowed with a topology $\tau$, such that the vector spaces operations, addition $+ \colon X \times X \to X$ and scalar multiplication $\cdot \co \K \times X \to X$, are continuous maps (see~\cites{Folland,Narici,Osborne,Rudin,Treves,Yosida}).
	\end{definition}

In a TVS $X$, the notion of Cauchy sequence (or Cauchy net, more generally) makes sense. A sequence $(x_{n})_{n\in \N}$ in $X$ is called a Cauchy sequence if the sequence $(x_{n}-x_{m})_{n,m\in \N}$ converges to zero in $X$. See Folland~\cite{Folland}.

A seminorm on a vector space $X$ is a function $p \co X \to \R$ such that:
	\begin{itemize}
		\item[(i)] $p(x)\geqslant 0$, for every $x\in X$.
		\item[(ii)] $p(\lambda x)=|\lambda| p(x)$, for every $x\in X$ and $\lambda \in \K$.
		\item[(iii)] $p(x+y)\leqslant p(x)+p(y)$, for every $x,y\in X$.
	\end{itemize}
The condition $(ii)$, in particular, says that $p(0)=0$ for all seminorm $p$. If $x=0$ whenever $p(x)=0$, we say that $p$ is a norm.

	\begin{definition}
		A separating family of seminorms on a vector space $X$ is a family of seminorms $(p_{\alpha})_{\alpha \in \Lambda}$ on $X$ such that for every nonnull vector $x \in X$ there exists $\alpha = \alpha_{x} \in \Lambda$ with $p_{\alpha}(x) \not = 0$.
	\end{definition}

If $(p_{\alpha})_{\alpha \in \Lambda}$ is a separating family of seminorms on a vector space $X$, they can be used to define a topology on $X$, in the same way we use the norm to define one. More precisely, given $x\in X$, $\alpha \in \Lambda$ and $\varepsilon >0$ we put
	\[
		B_{\alpha}(x;\varepsilon):=\{y\in X:p_{\alpha}(y-x)<\varepsilon  \}
	\]
so that the collection ${\mathcal B}$ of all finite intersection of sets $B_{\alpha}(x;\varepsilon)$ defines a base for a Hausdorff topology $\tau$ on $X$, which is called the \emph{topology generated by the family of seminorms} $(p_{\alpha})_{\alpha \in \Lambda}$. With this topology, a sequence $(x_{n})_{n\in \N}$ in $(X, \tau)$ converges to a vector $x\in X$ if and only if it converges to $x$ in every seminorm $p_{\alpha}$, that is, for every $\alpha\in \Lambda$ we have
	\[
		p_{\alpha}(x_{n}-x)\underset{n\to \infty}{\longrightarrow}0,
	\]
since the collection of sets of the form
	\[
		\bigcap_{j=1}^{n} B_{\alpha_{j}}(x;\varepsilon_{j})
	\] 
with $\alpha_{1},\cdots \alpha_{n}\in \Lambda$, $\varepsilon_{1},\cdots, \varepsilon_{n}>0$ and $n\in {\mathbb N}$, is a local base for $\tau$ at $x$.

In particular, a series $\sum x_{n}$ converges to $s\in X$ if and only if for every $\alpha \in \Lambda$
	\[
		p_{\alpha}\left(s-\sum_{j=1}^{n}x_{j}\right)\underset{n\to \infty}{\longrightarrow}0.
	\]
A series $\sum x_{n}$ is said to be \emph{absolutely convergent} if for every $\alpha \in \Lambda$ the series $\sum p_{\alpha}(x_{n})$ is convergent.

In~\cite{Babalola}, the author assumes that the family of seminorms are always \emph{saturated}, in the sense that for every finite collection $J$ of indexes $j$ we add $q_{J}:=\max_{j \in J} p_j$ to the family of seminorms. The topology generated by the saturated family coincides with the one generated by the original family of seminorms but has the advantage that a typical basic neighbourhood of $0$ is of the form $\{y\in X:q_{J}(y) <\varepsilon\}$, for some seminorm $q_{J}$; that is, no intersections needed. For our purpose, this will be necessary only on Theorem \ref{Th:generator}. We may change the family of seminorms $\{p_1, p_2, p_3, \ldots\}$ by $\{p_1, p_1+p_2, p_1+p_2+p_3, \ldots\}$ without altering the topology; in particular, $p_1 \leqslant p_1+p_2 \leqslant p_1 + p_2 + p_3 \leqslant \cdots$. See~\cites{Narici,Treves}.

	\begin{definition}
		A Fr\'{e}chet space $X= \big( X, (p_{j})_{j\in {\mathbb N}} \big)$ is a TVS $X$ whose topology is given by a countable separating family of seminorms $(p_{j})_{j\in {\mathbb N}}$ with the property that every Cauchy sequence converges in $X$.
	\end{definition}

Thus every Banach space is a Fr\'{e}chet space, although there are Fr\'{e}chet spaces which are not Banach spaces.

Note that if $\sum x_{n}$ is an absolutely convergent series in a Fr\'{e}chet space $X=\big(X,(p_{j})_{j\in {\mathbb N}}\big)$ then the sequence $\big( \sum_{k=1}^{n}x_{k} \big)_{n\in \N}$ is a Cauchy sequence in $X$ (consequently it is convergent), since for every $j\in \N$
	\[
		p_{j}\left(\sum_{k=m}^{n}x_{k}\right)\leqslant \sum_{k=m}^{n} p_{j}(x_{k})\underset{m,n\to \infty}{\longrightarrow}0.
	\]

	\begin{remark}
		Every Fr\'{e}chet space $X=\big( X, (p_{j})_{j\in {\mathbb N}} \big)$ is metrizable. More precisely, its topology is given by a translation-invariant metric $d\colon X\times X \to \R$, namely
			\[
				d(x,y):=\sum_{j=1}^{\infty}\frac{1}{2^{j}}\frac{p_{j}(x-y)}{1+p_{j}(x-y)}.
			\]
	\end{remark}

		As a metric space, the notion of boundness of a subset $B \subset X$ could be defined using the metric $d$, but we will not do this. With this metric, every subset would be bounded. We say that a subset $B$ of a Fr\'{e}chet space $X$ is \emph{bounded} if for every neighbourhood $V$ of the origin of $X$ there exists $t_{0}>0$ such that
	\[
		B\subset tV, \, \textrm{  for all } t \geqslant t_{0}.
	\]

Actually, this is the definition of boundedness on a general TVS $X$.

If the topology of $X$ is given by a separating family of seminorms $(p_{\alpha})_{\alpha \in \Lambda}$, this concept is equivalent to state that for every $\alpha \in \Lambda$ the set $\{p_{\alpha}(x):x\in B\}$ is bounded in ${\mathbb R}$, and we say that $B$ is bounded in every seminorm.

For now, let us see some examples of Fr\'{e}chet spaces. We abide by the convention that, unless otherwise stated, all functions are complex-valued.

	\begin{example} 
		Let $\Omega \subset {\mathbb R}^{n}$ be an open set, then $C^{\infty}(\Omega) = \big(C^{\infty}(\Omega), (p_{m,j})_{(m,j) \in \mathbb{Z}_{+} \times \mathbb{N}}\big)$ is a Fr\'{e}chet space, where
		\[
			p_{m,j}(\phi) := \sum_{|\alpha| \leqslant m} \sup_{K_j} \big| \partial^{\alpha}\phi \big|, \textrm{ for } \phi \in C^{\infty}(\Omega),
		\]and $K_j := \big\{x \in \mathbb{R}^{n}: |x| \leqslant j \textrm{ and } d(x, \partial\Omega) \geqslant 1/j \big\}$, so that $\{K_j\}_{j \in \mathbb{N}}$ exhausts $\Omega$. That is, the space of all infinitely differentiable functions endowed with the topology of the uniform convergence of the functions and their derivatives on compact subsets of $\Omega$ is a Fr\'{e}chet space,~\cite{Folland}.

		Let ${\mathscr E}'(\Omega)$ be the set of distributions with compact support on $\Omega$, that is, the topological dual space of $C^{\infty}(\Omega)$.
	\end{example}

	\begin{example}
		Let ${\mathscr S}(\R^{n})$ be the space of all functions $u\in C^{\infty}(\R^{n})$ such that
			\[
				\|u\|_{(N,\alpha)}:=\sup_{x\in \R ^{n}}(1+|x|)^{N} \big|\partial^{\alpha}u(x) \big|
			\]
		is finite for every non-negative integer $N$ and every multi-index $\alpha = (\alpha_{1},\cdots,\alpha_{n})\in \Z^{n}_{+}$. That is, ${\mathscr S}(\R^{n})$ is the space of all functions $u\in C^{\infty}(\R^{n})$ such that $x^{\alpha}\partial ^{\beta}u$ goes to zero when $|x|\to \infty$, for all multi-indexes $\alpha, \beta$.
		
		It follows that $\big( \|\cdot \|_{(N,\alpha)} \big)_{N,\alpha}$ is a countable separating family of seminorms on ${\mathscr S}(\R^{n})$ which turns it into a Fr\'{e}chet space (see Folland~\cite{Folland}).
		
		More generally, given Hilbert spaces $H_0 = (H_0, \| \cdot \|_{H_0})$ and $H_1 = (H_1, \| \cdot \|_{H_1})$, we may consider functions $\phi \colon \Omega \subset \R^n \to (H_0,\| \cdot \|) $ and then consider the Fr\'{e}chet spaces $C^{\infty}(\Omega, H_0) = \big(C^{\infty}(\Omega, H_0), (p_{m,j})_{(m,j) \in \mathbb{Z}_{+} \times \mathbb{N}}\big)$ and ${\mathscr S}(\R^{n},H_0) = \big( {\mathscr S}(\R^{n},H_0), (\|\cdot \|_{(N,\alpha)} )_{N,\alpha} \big)$, where the seminorms are adapted obviously considering the norm of $H_0$.
		
		Let ${\mathscr S}'(\R^n, H_0; H_1)$, or ${\mathscr S}'$ for short, denote the space of all bounded linear applications from ${\mathscr S}(\R^{n},H_0)$ into $H_1$; similarly, ${\mathscr E}'(\Omega,H_0;H_1)$ is defined. If $H_0=H_1=\C$ then we get the usual spaces of distributions ${\mathscr S}'(\R^{n})$ and ${\mathscr E}'(\Omega)$.

	\end{example}

The following example is less known,~\cites{aragao1,Treves1976}, and it was generalized in the sense described above. It is noteworthy that the Fourier transform ${\mathscr F}$ is well defined from ${\mathscr S}'(\R^n, H_0; H_1)$ into ${\mathscr S}'(\R^n, H_0; H_1)$ by setting
	\[
		({\mathscr F}u)\phi = \langle {\mathscr F}u, \phi \rangle := \langle u, {\mathscr F}\phi \rangle \in H_1,
	\]
for every $u \in {\mathscr S}'(\R^n, H_0; H_1)$ and $\phi \in {\mathscr S}(\R^{n},H_0)$. We shall frequently write $\widehat{u}$ instead of ${\mathscr F}u$.

	\begin{example}\label{FL2loc}
		Let $L^{2}_{loc}\big(\R^{n}, {\mathscr L}(H_0,H_1) \big)$, or $L^{2}_{loc}$ for short, be the set of all measurable functions with values on ${\mathscr L}(H_0,H_1)$ whose norm square is integrable on every compact set of $\R^{n}$. For simplicity, the reader may assume that $H_0=H_1=\C$. Set
			\[
				E := \big\{ u \in {\mathscr S}'(\R^{n},H_0;H_1): \widehat{u} \in L^{2}_{loc}\big(\R^{n}, {\mathscr L}(H_0,H_1) \big) \big\}
			\]
		and endow it with the topology generated by the seminorms 
		\begin{equation}\label{seminorm-E}
			p_{j}^{*}(u) := \left (\int_{|\xi|\leqslant j} \| \widehat{u}(\xi) \|^{2}_{{\mathscr L}(H_0,H_1)} \, d\xi \right)^{1/2}, \textrm{ for } u\in E \textrm{ and } j\in \N.
		\end{equation}
		
		It follows that this family of seminorms on $E$ is a separating one, whence the function
			\[
				d(u,v) := \sum_{j=1}^{\infty}\frac{1}{2^{j}}\frac{p_{j}^{*}(u-v)}{1+p_{j}^{*}(u-v)}
			\]
		defines a metric on it.
		
		Let ${\mathscr F}L^{2}_{loc}(\R^{n}, H_0; H_1)$, or ${\mathscr F}L^{2}_{loc}$ for short, be the completion of the metric space $(E,d)$. If $[u]\in {\mathscr F}L^{2}_{loc}$, we can define its Fourier transform\footnote{Here, we are considering the completion as the quotient space of all Cauchy sequence in $E$ under the canonical equivalence: $(u_{l})_{l\in \N} \sim (v_{l})_{l\in \N} \Longleftrightarrow \underset{l\to \infty}{\lim}d(u_{l},v_{l})=0.$}: if $(u_{l})_{l\in \N}\in [u]$ then $(\widehat{u_{l}})_{l\in \N}$ is a Cauchy sequence in $L^{2}_{loc}$, consequently there exists a unique $w \in L^{2}_{loc}$ such that $\widehat{u_{l}} \overset{l\to \infty}{\longrightarrow}w$ in $L^{2}_{loc}$ and we set $\widehat{[u]} := w$, so that the Fourier transform maps ${\mathscr F}L^{2}_{loc}$ into $L^{2}_{loc}$:
		\[
			{\mathscr F} \colon {\mathscr F}L^{2}_{loc} \to L^{2}_{loc}.
		\]
		
		It is straightforward to see that $\widehat{[u]}$ is independent of the sequence $(u_{l})_{l\in \N}$ chosen to represent the class $[u]$, so the Fourier transform is well defined. Maybe we should have denoted this space by ${\mathscr F}^{-1}L^{2}_{loc}$ so the Fourier transform above would have been appropriately defined, in terms of notation. But we follow the original notation, found in~\cite{Treves1976}.
		
		If $p_{j}^{*} \colon {\mathscr F}L^{2}_{loc} \to [0, \infty)$ denotes the natural extensions of the seminorms $(\ref{seminorm-E})$ to ${\mathscr F}L^{2}_{loc}$, it is not hard to see that topology of the complete (metric, and also linear) space ${\mathscr F}L^{2}_{loc}$ is equivalent to the topology generated by $(p_{j}^{*})_{j \in \N}$ and hence ${\mathscr F}L^{2}_{loc} = \big( {\mathscr F}L^{2}_{loc}, (p_{j}^{*})_{j \in \N} \big)$ is a Fr\'{e}chet space,~\cite{Treves1976}.
	\end{example}

{\color{red}
\begin{remark}\label{Obs:FL2loc}
It is quite important to note that $L^{2}=L ^{2}\big(\R^{n}, {\mathscr L}(H_0,H_1) \big)$ and ${\mathscr E}'\big(\R^{n}, {\mathscr L}(H_0,H_1) \big)$ are natural subspaces of ${\mathscr F}L^{2}_{loc}\big(\R^{n}, {\mathscr L}(H_0,H_1) \big)$.

Indeed, $L^2 \subset {\mathscr F}L^{2}_{loc}$ simply because the Fourier transform is an isometric isomorphism from $L^{2}$ onto itself, by Plancherel Theorem~\cites{Folland,Folland1992,Stein,Treves}. Moreover, by Paley-Wiener-Schwartz Theorem~\cites{Stein,Treves}, a temperated distribution $v$ on $\R^{n}$ belongs to ${\mathscr E}'$ if and only if $\widehat{v}$ has an entire extension $V:\C^{n} \to {\mathscr L}(H_0,H_1)$. In particular, ${\mathscr E}'$ is also a subspace of ${\mathscr F}L^{2}_{loc}$.
\end{remark}
}

Although it is quite general to consider ${\mathscr F}L^{2}_{loc}(\R^{n}, H_0; H_1)$, there is no loss of generality by assuming that $H_0 =H_1 =\C$, so we do it from now on. The reader is invited to verify that every proof (concerning the Fr\'{e}chet space ${\mathscr F}L^{2}_{loc} = {\mathscr F}L^{2}_{loc}(\R^n,\C; \C)$) may be naturally adapted to the general case.

\begin{definition} Let $T\colon X \to X$ be a linear operator on $X$. We say that $T$ is a bounded linear operator if it takes bounded sets of $X$ to bounded set of $X$.

We denote the space of all bounded linear operators from $X$ into $X$ by ${\mathscr L}(X)$.
\end{definition}

If $T\colon X\to X$ is a linear operator then the following statements are equivalent (you may see~\cites{Folland,Narici,Osborne,Rudin,Treves,Yosida} for a proof): 
\begin{itemize}
\item[(i)] $T\in {\mathscr L}(X)$.
\item[(ii)] $T$ is continuous.
\item[(iii)] For every $j\in {\mathbb N}$ there exist indexes $l_{1}, \cdots, l_{k}\in {\mathbb N}$ and a constant $C>0$ (all depending only on $j$ and $T$) such that
	\[
		p_{j}(Tx)\leqslant C \sum_{r=1}^{k}p_{l_{r}}(x), \textrm{ for all } x \in X.
	\]
\end{itemize}

We restrict our attention to those pseudodifferential operators (\cites{HormanderIII,FollandPDE}) with symbols independent of the space variable $x$ defined on ${\mathscr S}(\R^{n})$ in order to obtain the groundwork over which the theory of generation of groups can be applied to solve evolution problems.

\begin{definition}\label{Def.:pseudodiff}
A pseudodifferential operator of order $m$ on $\mathbb{R}^n$ with constant coefficients (or constant coefficients $m$-$\Psi$DO for short) is a linear map $a(D) \colon {\mathscr S}(\R^{n}) \to {\mathscr S}(\R^{n})$ given by
	\[
		\big( a(D)\psi \big)(x) := \int_{\mathbb{R}^n} e^{2 \pi i x \cdot \xi}a(\xi) \widehat{\psi}(\xi) \, d\xi, \textrm{ for every } x \in \mathbb{R}^{n},
	\]
where $a \in C^{\infty}(\mathbb{R}^n)$ satisfies the property that for all multiindex $\alpha$ there is a constant $c_{\alpha}>0$ such that
	\[
		\big| \partial^{\alpha} a(\xi) \big| \, \leqslant \, c_{\alpha} (1+ |\xi|)^{m-|\alpha|}, \, \xi \in \mathbb{R}^n.
	\]
\end{definition}

We say that $a$ is a symbol of order $m$ on $\mathbb{R}^n$ and the space of those functions, denoted by $S^m(\mathbb{R}^n)$, is a Fr\'{e}chet space with seminorms given by the smallest constants $c_{\alpha}$ which can be used in the inequality above. This definition also yields a one-to-one correspondence between constant coefficients $m$-$\Psi$DOs and symbols $a \in S^m(\mathbb{R}^n)$.

By the Proposition $8.3$ and Corollary $8.23$ of Folland~\cite{Folland} and by the Leibniz formula, every constant coefficients $m$-$\Psi$DO $a(D)$ is a continuous linear operator from ${\mathscr S}(\R^{n})$ to itself.

As the reader may readily see, every constant coefficients $m$-$\Psi$DO $a(D)$ induces an operator on the space of tempered distributions ${\mathscr S}'(\R^n)$ by setting
	\[
		\langle a(D)u, \psi \rangle := \langle u, a(D)\psi \rangle,
	\]
for every $u \in {\mathscr S}'(\R^n)$ and $\psi \in {\mathscr S}(\R^n)$.

Throughout this paper  $\hat{u}(\xi) = \int_{\mathbb{R}^n} e^{-2\pi i x \cdot \xi}u(x) \, dx$ shall denote the Fourier transform of $u \in L^{1}(\mathbb{R}^n)$ and we adopt the notational convention that $D = \frac{1}{2\pi i}\partial$, so that the convenient formula $(D^{\alpha}u) \hat{ \ }(\xi) = \xi^{\alpha}\hat{u}(\xi)$ holds, for all multiindex $\alpha$.

The class of constant coefficients $m$-$\Psi$DOs is a large one which contains every linear differential operator
	\[
		a(D) = \sum_{|\alpha|\leqslant m} a_{\alpha}D^{\alpha} \colon {\mathscr S}(\R^{n}) \to {\mathscr S}(\R^{n}),
	\]
where $a_{\alpha} \in \mathbb{C}$. To see this, just apply the Fourier inversion theorem to the formula $(D^{\alpha}u)\hat{ \, }(\xi) = \xi^{\alpha}\hat{u}(\xi)$, for $u \in {\mathscr S}(\R^{n})$, to get
\[
	\big( a(D)u \big)(x) = \sum_{|\alpha|\leqslant m} a_{\alpha} \int_{\mathbb{R}^n} e^{2 \pi i x \cdot \xi} \xi^{\alpha} \hat{u}(\xi) \, d\xi
								 = \int_{\mathbb{R}^n} e^{2 \pi i x \cdot \xi} a(\xi) \hat{u}(\xi) \, d\xi,
\]
where $a(\xi) = \sum a_{\alpha} \xi^{\alpha}$ is the symbol of $a(D)$. More generally, it is straightforward to check that the class of pseudodifferential operators associated to the symbols $\mathbb{R}^n \times \mathbb{R}^n \ni (x,\xi) \mapsto a(x,\xi) \in \C$ contains all linear differential operators $a(x,D) = \sum_{|\alpha|\leqslant m} a_{\alpha}(x) D^{\alpha} \colon {\mathscr S}(\R^{n}) \to {\mathscr S}(\R^{n})$, where $a_{\alpha} \in C^{\infty}(\mathbb{R}^n)$. See~\cites{FollandPDE,HormanderIII}.


\section{\bf{Generation theorem and consequences}}\label{Section:generation}

At this point, the strategy for solving the problem \eqref{eq:Cauchy} in a Fr\'{e}chet space $X$ is pretty clear. The addition on $X$ and scalar multiplication are well defined operations and are continuous with the obvious product topologies. We need to require a continuous linear operator $A$ defined from this space to itself to have the appropriate compatibility with the topology on $X$ so its exponential operator $\exp(A)\colon X \to X$ makes sense and may be used to define the solution of the associated Cauchy problem. Such strong compatibility is expressed in terms of the seminorms on the space in Definition~\ref{Def:LscX} and fortunately constant coefficients $m$-$\Psi$DOs defined on ${\mathscr F}L^{2}_{loc}$ will naturally satisfy theses conditions, as we shall see in Theorem~\ref{Th:generationFL2loc}.

Let us make it precise. Fix a Fr\'{e}chet space $X = \big( X, (p_{j})_{j \in \mathbb{N}} \big)$.

	\begin{definition}\label{Def:LscX}
		A bounded linear operator $A\in {\mathscr L}(X)$ is said to be strongly compatible with $(p_{j})_{j \in \mathbb{N}}$ if satisfies the following properties:
			\begin{enumerate}
				\item[i)] $p_{j}(Ax)=0$ whenever $p_{j}(x)=0$, for every $j\in \N$; and
				\item[ii)] $\sup \big\{p_{j}(Ax): p_{j}(x)=1 \big\}$ is finite, for every $j \in \N$.
		\end{enumerate}
		
		We denote by $\Lsc$ the set of all $A\in {\mathscr L}(X)$ which are strongly compatible with $(p_{j})_{j \in \mathbb{N}}$.
	\end{definition}

Note that if $p\colon X\to \R$ is a norm then every $T\in {\mathscr L}(X)$ is strongly compatible with $p$. It is noteworthy that $\Lsc$ is not empty since at least the identity operator of $X$ is strongly compatible with $(p_{j})_{j \in \mathbb{N}}$.

About the requirements above, the first one has appeared in~\cite{aragao1} as a natural condition to obtain exponential dichotomy for evolution processes and it certainly does not turn $p_j$ into a norm; it is a weaker requirement over the operator, not over the seminorms. An operator $A\in {\mathscr L}(X)$ that satisfies i. is said to be compatible with $(p_{j})_{j \in \mathbb{N}}$.

It is a simple exercise to verify that, if $A \in {\mathscr L}(X)$ is compatible with $(p_{j})_{j \in \mathbb{N}}$ then one of the expressions below
	\begin{equation}\label{eq:supremum}
		\sup_{p_{j}(x)=1} p_{j}(Ax), \sup_{p_{j}(x)<1} p_{j}(Ax) \mbox{ and } \sup_{p_{j}(x) \leqslant 1} p_{j}(Ax)
	\end{equation}
is finite if and only if all three are; and in this case, they all coincide. In particular, thanks to Lemma \ref{Lemma:seminorms}, there exists a positive constant $c=c(j,A)$ such that $A\big( B_j(0,1) \big) \subset c B_j(0,1)$; according to the notation in Section \ref{Section:preliminaries}.

The second condition on the Definition \ref{Def:LscX} is not obvious, since $\{x \in X: p_j(x)=1\}$ is not in general a bounded set on $X$, consequently $A\big(\{x \in X: p_j(x)=1\}\big)$ need not to be bounded as well. Also, it is needed to topologyze appropriately $\Lsc$ in order to define the operator $\exp(A)$ in $\Lsc$. We define seminorms on $\Lsc$ by setting \footnote{Obviously, a seminorms $p$ such that $p(x) = 0$ for all $x \in X$ is a null seminorm, so that we have no interest on such $p$ and hence it will be dismissed from the family of seminorms on $X$.}
\[
	p_{j}^{X}(A):=\sup_{p_{j}(x)=1} p_{j}(Ax), \, j \in \N,
\]
for every $A \in \Lsc$.

\begin{remark}\label{Remark:Babalola-condition}
The requirements in Definition \ref{Def:LscX} are closely related to the one that has already appeared in~\cite{Babalola}. In fact, they are equivalent. This is the key point in common with Babalola's approach, but as we shall see in Remark~\ref{Remark:comparing-to-Babalola} it is practically the only one.

Let $Y$ be a complete Hausdorff locally convex TVS over $\C$ and let $\{q_{\lambda}\}_{\lambda \in \Lambda}$ be a saturated family of seminorms on it. Set $V_{\lambda} := \{y \in Y: q_{\lambda}(y) < 1\}$ for every $\lambda \in \Lambda$. The author deals with a class of linear operators $S \colon Y \to Y$ with the property that, for every $\lambda \in \Lambda$, there exists a positive constant $c=c(\lambda,S)$ such that
	\[
		S V_{\lambda} \subset c V_{\lambda},
	\]
and he writes $S \in {\mathscr L}_{A}(Y)$. Hence, $S$ is bounded and satisfies $q_{\lambda}(Sy) \leqslant c(\lambda,S) q_{\lambda}(y)$ for every $y \in Y$ and every $\lambda \in \Lambda$.

In particular, if $Y=\big(Y, (q_j)_{j \in \N} \big)$ is a Fr\'{e}chet space, with $q_j \leqslant q_{j+1}$ for every $j$, and $S \in {\mathscr L}_{A}(Y)$ then $S$ is a strongly compatible operator; that is, the Babalola's condition implies ours. On the other hand, by \eqref{eq:supremum}, every strongly compatible operator $S\colon Y \to Y$ satisfy the Babalola's condition.
\end{remark}

One can easily verify that $\Big( \Lsc, \big( p_{j}^{X} \big)_{j \in \mathbb{N}} \Big)$ is a Fr\'{e}chet space, thanks to the following lemma.


\begin{lemma}\label{Lemma:seminorms}
If $A\in \Lsc$ then $p_{j}(Ax)\leqslant p_{j}^{X}(A)p_{j}(x)$, for every $x \in X$ and every $j$.

In particular, for every $n \in \mathbb{N}$, there holds
	\[
		p_{j}^{X}(A^{n})\leqslant p_{j}^{X}(A)^{n}.
	\]
\end{lemma}
\proof
Indeed, for a fixed $j\in \N$, the inequality $p_{j}(Ax) \leqslant p_{j}^{X}(A)p_{j}(x)$ is trivial whenever $p_{j}(Ax)=0$.

If $x \in X$ satisfies $p_{j}(Ax) \neq 0$ then $p_{j}(x) \neq 0$ and for $x_0 = \dfrac{1}{p_{j}(x)} \, x$ we get
\[
	\dfrac{p_j(Ax)}{p_j(x)} = p_{j}(Ax_0) \leqslant \sup_{p_{j}(z) =1} p_{j}(Az)=p_{j}^{X}(A)
\]
whence $p_{j}(Ax) \leqslant p_{j}^{X}(A) \, p_j(x)$, for every $x \in X$.

Furthermore, since $p_{j}(A^n x) \leqslant p_{j}^{X}(A) \, p_{j}(A^{n-1}x)$, for every $x \in X$ and every natural number $n \geqslant 2$, the result then follows by induction.

\endproof

If $Y=\big(Y, (q_k)_{k \in \N} \big)$ is another Fr\'{e}chet space then we could have defined the space ${\mathscr L}_{\textrm{sc}}(X,Y)$ and proved the previous lemma to it, with the obvious adaptations.

Unsurprisingly, we extend the concepts of group of bounded linear operators, $C_0$-group and uniformly continuous group on Fr\'{e}chet spaces.

\begin{definition}
A family $\{T(t):t \in \R\} \subset {\mathscr L}(X)$ is called a group of bounded linear operators on $X$ (or a group on $X$, for short) if:
\begin{itemize}
\item[(i)] $T(0)=I$, where $I\colon X\to X$ is the identity operator on $X$.
\item[(ii)] $T(t+s)=T(t)T(s)$, for every $t,s \in \R$.
\end{itemize}

We say that a group $\{T(t): t \in \R\} \subset {\mathscr L}(X)$ is a $C_{0}$-group if 
\begin{itemize}
\item[(iii)] $T(t)x \overset{X}{\underset{t \to 0}{\longrightarrow}} x$, for every $x\in X$.
\end{itemize}

A group $\{T(t):t \in \R\} \subset \Lsc$ is said to be a uniformly continuous group if
\begin{itemize}
\item[(iv)] $T(t)\overset{\Lsc}{\underset{t\to 0}{\longrightarrow}} I$, that is,  $p_{j}^{X}\big(T(t)- I\big) \overset{\R}{\underset{t\to 0}{\longrightarrow}} 0$ for every $j \in \N$.
\end{itemize}
\end{definition}

\begin{definition}
Let $\{T(t):t \in \R\} \subset {\mathscr L}(X)$ be a semigroup on $X$. Its infinitesimal generator is the linear operator $A\colon D(A)\subset X \to X$ defined by 
\[
	Ax := \lim_{t \to 0} \dfrac{T(t)x-x}{t},
\]
for $x \in D(A) := \left\{ x \in X: \textrm{ the limit } \underset{t\to 0}{\lim}\frac{T(t)x-x}{t} \textrm{ exists in } X \right\}$.
\end{definition}

We often will write $T(\cdot)$ to denote the group $\{T(t): t \in \R\}$ on $X$. Clearly, a group $T(\cdot)$ has a unique infinitesimal generator. If we replace $\R$ by the interval $[0,\infty)$ in the last two definitions, we get the concept of semigroups of bounded linear operator on $X$, indicated by $\{ T(t):t\geqslant 0 \}$, and its infinitesimal generator $A$ whose definition of course is obtained replacing $\displaystyle\lim_{ t\to 0}$ by $\displaystyle\lim_{t\to 0^{+}}$.

It is immediate to see that these definitions coincide with the usual ones if $X$ is a Banach space.

Next we state and prove the main results of this paper concerning the generation of uniformly continuous groups on $X$. By Pazy~\cite{Pazy}, Theorem $1.2$, we know that, provided that $X$ is a Banach space, a linear operator $A\colon D(A) \subset X \to X$ is the infinitesimal generator of a uniformly continuous group if and only if $A$ is a bounded linear operator on $X$. The sufficiency is quite immediate as we state and prove below; and it is a consequence of the previous lemma. We will refer to it as simply the generation theorem.

\begin{theorem}\label{Th:generation} For every $A\in \Lsc$ and $t\in \R$ the series
\[
	\sum_{n=0}^{\infty}\frac{(tA)^{n}}{n!}
\]
converges in $\Lsc$ and, naturally, we indicate its sum by $e^{tA}.$

Besides, the family $\{ e^{tA}:t\in \R \}$ is a uniformly continuous group of bounded linear operators on $X$ and the operator $A$ is its infinitesimal generator.
\end{theorem}
\proof For each $N\in \N$ let
\[
	S_N := \sum_{n=0}^{N} \dfrac{(tA)^n}{n!} \in \Lsc,
\]
and let us show that for each fixed $j\in \N$ the sequence $(S_{N})_{N\in \N}$ is a Cauchy sequence with respect to the seminorm $p_{j}^{X}$.

Indeed, given $\varepsilon>0$, if $N>M$ are natural numbers we have
\[
	p_{j}^{X}(S_N - S_M) = p_{j}^{X}\left( \sum_{n=M+1}^{N} \dfrac{(tA)^n}{n!} \right)
									   \leqslant \sum_{n=M+1}^{N} \dfrac{\left( t \, p_{j}^{X}(A)\right)^n}{n!}
									   < \varepsilon, 
\]
for  $N,M$  large enough, since the numerical series $\sum \frac{t^{n}\left( \, p_{j}^{X}(A)\right)^n}{n!}$ converges to the $e^{t p^{X}_{j}(A)}$ in $\R$.

Now, it is clear that $e^{0A}$ is the identity of $X$.

Also, since $\sum_{n=0}^{\infty}\frac{(tA)^{n}}{n!}$ is absolutely convergent by the classical formula for product of series (see Bartle~\cite{Bartle}, Theorem $26.15$), we conclude that $e^{(s+t)A} = e^{sA}e^{tA}$ for all $t,s \in \R$. Indeed,
\begin{align*}
	e^{(s+t)A} &= \sum_{n=0}^{\infty}\frac{[(s+t)A]^{n}}{n!}
				  = \sum_{n=0}^{\infty} \left[\sum_{k=0}^{n} \left( \begin{array}{c} n \\ k \end{array} \right) s^{n-k}t^{k} \right] \frac{A^{n}}{n!}\\
		&= \sum_{n=0}^{\infty} \sum_{k=0}^{n} \frac{(sA)^{n-k}}{(n-k)!} \circ \frac{(tA)^{k}}{k!}
		=\left(\sum_{n=0}^{\infty}\frac{(sA)^{n}}{n!}\right)\circ \left(\sum_{n=0}^{\infty}\frac{(tA)^{n}}{n!}\right)
		= e^{sA}e^{tA}.
\end{align*}

Besides, for every $j\in \N$ we have
\[
	p^{X}_{j}(e^{tA}-I) = p^{X}_{j}\left( \sum_{n=1}^{\infty} \dfrac{(tA)^n}{n!} \right)
									\leqslant \sum_{n=1}^{\infty}\frac{\big(t \, p^{X}_{j}(A)\big)^{n}}{n!}
									= e^{t p^{X}_{j}(A)}-1,
\]
whence $\{ e^{tA}:t\in \R \}$ is a uniformly continuous group of bounded linear operators on $X$.

Finally, by the definition of generator, fixed $j\in \N$, if $x\in X$ and $t\not = 0$ we have
\[
	p_{j}\left(\frac{e^{tA}x-x}{t} -Ax\right)
		\leqslant \frac{1}{t}\sum_{n=2}^{\infty}\frac{\big(t \, p^{X}_{j}(A)\big)^{n}}{n!}p_{j}(x)
		= \left( \frac{e^{t p^{X}_{j}(A)} -1 }{t} -p^{X}_{j}(A) \right) p_{j}(x),
\]
hence $A$ is the infinitesimal generator of $\{e^{tA}:t\in \R\}$.

\endproof

If $X$ is a Banach space then every bounded linear operator $A \colon X \to X$ is the infinitesimal generator of a unique uniformly continuous group on $X$. Is the same true for Fr\'{e}chet spaces? The answer is affirmative and it is a straightforward generalization of the proof for Banach spaces.

\begin{proposition}\label{unico-semigrupo}
If $T(\cdot)$ and $S(\cdot)$ are uniformly continuous groups on $X$ so that
	\[
		\lim_{t \to 0} \dfrac{T(t)-I_X}{t} = A = \lim_{t \to 0} \dfrac{S(t)-I_X}{t} \mbox{ in } \Lsc,
	\]
then $T(t)=S(t)$ for every $t \in \R$.
\end{proposition}

\proof

We just have to prove that, given $\tau>0$, $T(t) = S(t)$ for every $0 \leqslant t \leqslant \tau$; or equivalently, given $\tau>0$ and given $\varepsilon>0$, we get $p_j^{X}\big( T(t) - S(t) \big) \leqslant \varepsilon$ for every $0 \leqslant t \leqslant \tau$ and $j \in \N$.

By continuity of the maps $t \mapsto p_j^{X}\big( S(t) \big)$ and $t \mapsto p_j^{X}\big( T(t) \big)$, there exists a positive constant $c=c(j,\tau, S, T)>0$ such that
	\[
		\sup_{0 \leqslant s,t \leqslant \tau} p_j^{X}\big( T(t) \big) \, p_j^{X}\big( S(s) \big) \leqslant c.
	\]

By hypothesis, there exists a positive constant $\delta = \delta(j, \tau, \varepsilon, S, T)>0$ such that
	\[
		\sup_{0 \leqslant h \leqslant \delta} h^{-1} p_j^{X} \big( T(h) - S(h) \big) \leqslant \dfrac{\varepsilon}{\tau c}.
	\]

For $0\leqslant t \leqslant \tau$, take $n \in \N$ such that $t/n < \delta$, so
	\begin{align*}
		p_j^{X}\big( T(t) - S(t) \big) & = p_j^{X}\left( T\left(n \frac{t}{n} \right) - S\left(n \frac{t}{n}\right) \right) \\
			& \leqslant \sum_{k=0}^{n-1} p_j^{X} \left(
					T \left( \dfrac{(n-k)t}{n} \right) S\left( \dfrac{kt}{n} \right)
					- T \left( \dfrac{(n-k-1)t}{n} \right) S\left( \dfrac{(k+1)t}{n} \right)
				\right) \\
			& \leqslant \sum_{k=0}^{n-1}
				p_j^{X} \left( T \left( \dfrac{(n-k-1)t}{n} \right) \right)
				p_j^{X} \left( T\left(\frac{t}{n} \right) - S\left(\frac{t}{n}\right) \right)
				p_j^{X} \left( S\left( \dfrac{kt}{n} \right) \right) \\
			& \leqslant n \, c \, \dfrac{\varepsilon}{\tau c} \dfrac{t}{n} \leqslant \varepsilon
	\end{align*}
and the proof is complete.

\endproof

Now we stablish some kind of reciprocal of the generation theorem. To do so, we must fix some notation: for every $j \in \N$, set $X_j := \big( X/p_j^{-1}\big(\{0\}\big), \| \cdot \|_j)$, where
	\[
		\| [x]_j \|_j := \inf_{p_j(z)=0} p_j(x-z),
	\]
for $[x]_j$\footnote{Of course, $[x]_j$ stands for the equivalence class of $x\in X$ relatively to the quotient space $X/p_j^{-1}\big(\{0\}\big)$.} in $X/p_j^{-1}\big(\{0\}\big)$, whence every $X_j$ is a normed space.

\begin{theorem}\label{Th:generator}
If $\{T(t): t \in \R\}\subset \Lsc$ is a uniformly continuous group on $X$ and every $X_j$ is a Banach space, then its infinitesimal generator, namely $A\colon D(A) \subset X \to X$, is defined on whole space $X$, $A \in \Lsc$ and $T(t) = e^{tA}$ in $\Lsc$, for every $t \in \R$.
\end{theorem}

\proof
We may assume that the seminorms of $X$ are a nested increasing sequence, that is, $p_j \leqslant p_{j+1}$ for every $j$, as discussed earlier; hence $X_1 \subset X_2 \subset \cdots \subset X$.

First we induce a family of linear operators on every $X_j$ by setting
\begin{align*}
	T_j(t) \colon & X_j \to X_j \\
	&[x]_j \mapsto [T(t)x]_j
\end{align*}
for every $t \in \R$.

\medskip
\textbf{\underline{Claim 1:}} $\{T_j(t): t \in \R\}$ is a uniformly continuous group on $X_j$, for every $j$.
\medskip

Indeed, for fixed $j \in \N$, we have
	\begin{align*}
		\| T_j(t)[x]_j \|_j &= \inf_{p_j(z)=0} p_j\Big( T(t)x - T(t)z - \big(z-T(t)z\big) \Big) \\
										& \leqslant \inf_{p_j(z)=0} \left \{p_j^{X}\big(T(t)\big) p_j(x-z) + p_j(z) + p_j\big(T(t)z\big) \right \}\\
										&= p_j^{X}\big(T(t)\big) \| [x]_j \|_j,
	\end{align*}
since $T(t)$ is strongly compatible with $(p_k)_{k \in \N}$. So $T_j(t)$ is a bounded linear operator.

It is clear that $T_j(0)$ is the identity operator on $X_j$. Also, for $t,s \in \R$, we get
	\[
		T_j(t)\big(T_j(s)[x]_j\big) = [T(t) \circ T(s)x ]_j = T_j(t+s)[x]_j
	\]
and
	\begin{align*}
		\| T_j(t) - I_{X_j} \|_{\mathscr{L}(X_j)} &= \sup_{ \| [x]_j \|_j =1} \| T_j(t)[x]_j - [x]_j \|_j \\
						&= \sup_{\| [x]_j \|_j =1} \inf_{p_j(z)=0} p_j\big( T(t)x - x - z \big)\\
						&\leqslant \sup_{\| [x]_j \|_j =1} \inf_{p_j(z)=0} p_j\Big( \big(T(t)-I_X \big)(x - z) \Big) + p_j\big( T(t)z -z -z \big)\\
						&\leqslant \sup_{\| [x]_j \|_j =1} \inf_{p_j(z)=0} p_j^{X}\big(T(t)-I_X \big)p_j(x-z) \\
						&= p_j^{X}\big(T(t)-I_X\big) \ {\overset{\R}{\underset{t \to 0}{\longrightarrow}}} \ 0.
	\end{align*}

By Pazy~\cite{Pazy}, we know that
	\[
		T_j(t) = e^{tA_j} = \sum_{n=0}^{\infty} \dfrac{t^n}{n!} A_j^n,
	\]
with convergence in ${\mathscr L}(X_{j})$, where
	\[
		A_j := \big( T_j(t_j) - I_{X_j} \big) \circ \left( \int_0^{t_j} T_j(t) \, dt \right)^{-1} \in \mathscr{L}(X_j)
	\]
and $t_j>0$ is chosen so that $\|T_j(t) - I_{X_j} \|_{\mathscr{L}(X_j)} \leqslant 1/2$ for $0 \leqslant t \leqslant t_j$. It is noteworthy that the definition of $A_j$ does not depend on $t_j$ thanks to the uniqueness of the infinitesimal generator.

Let $\sigma_j \colon X \to X_j$ be the (continuous) canonical projection, that is, $\sigma(x) := [x]_j$, for every $j$; and set $\pi_j \colon X_{j+1} \to X_j$ by $\pi_j\big( [x]_{j+1} \big): = [x]_j$, which is well defined since $p_j \leqslant p_{j+1}$; and it is continuous, since for every $w \in X$ such that $p_{j+1}(w) =0$ we have
	\[
		\| \pi_j ([x]_{j+1}) \|_j = \inf_{p_j(z)=0} p_j(x-z) \leqslant \inf_{p_j(z)=0} \left\{p_j(x-w) + p_j(z-w) \right \}= p_j(x-w) 
	\]
and consequently 
	\[
		\| \pi_j ([x]_{j+1}) \|_j \leqslant \inf_{p_{j+1}(w)=0} p_{j+1}(x-w) =  \| [x]_{j+1} \|_{j+1}.
	\]

By construction, we get $\big( T_j(t) \circ \pi_j \big)([x]_{j+1}) = \big( \pi_j \circ T_{j+1}(t) \big)([x]_{j+1})$ and the diagram
	\begin{equation}\label{diagram:D1}
		\begin{split}
			\xymatrix{
			& X \ar[dr]^{\sigma_{j+1}} \ar[dl]_{\sigma_j} & \\
			X_j \ar[dd]_{A_j} & & X_{j+1} \ar[dd]^{A_{j+1}} \ar[ll]^{\pi_j} \\
			& & \\
			X_j & &  X_{j+1} \ar[ll]^{\pi_j} }	
		\end{split}
	\end{equation}
is commutative, for every $j \in \N$.

It is natural to seek the infinitesimal generator of $T( \cdot )$ using the infinitesimal generators $A_j$ of $T_j(\cdot)$, wondering whether exists a linear operator $A\colon X \to X$ such that every $A_j \colon X_j \to X_j$ is just the projection of $A$ on $X_j$ induced by $\sigma_j$; that is, $[Ax]_j = A_j[x]_j$ holds for every $j$ and $x\in X$. Well, this is the case.

\medskip
\textbf{\underline{Claim 2:}} there exists a unique linear operator $A\colon X \to X$ such that the diagram
	\begin{equation}\label{diagram:D2}
		\begin{split}
			\xymatrix{
			X \ar[r]^{A} \ar[d]_{\sigma_j} & X \ar[d]^{\sigma_j} \\
			X_j \ar[r]_{A_j}  &  X_j}
		\end{split}
	\end{equation}
is commutative, for every $j \in \N$.
\medskip

Fix $x \in X$.

Since every $\sigma_j$ is surjective, we obtain a sequence $(z_j)_j$ in $X$ such that $\sigma_j(z_j) = A_j \circ \sigma_j(x)$ for every $j\in \N$, and then
	\[
		\sigma_j(z_j) = A_j \circ \sigma_j(x)
						   = \pi_j\big(A_{j+1} \circ \sigma_{j+1}(x) \big)
						   = \pi_j \big( \sigma_{j+1}(z_{j+1}) \big)
						   = \sigma_j(z_{j+1})
	\]
so $\sigma_j(z_j - z_{j+1}) = [0]_j$, for every $j$; that is, $p_j(z_j-z_{j+1})=0$, for every $j$.

Since $p_j \leqslant p_{j+1}$, we get $p_{l}(z_j - z_k) = 0$ whenever $j,k \geqslant l$, that is, $(z_j)_j$ is a Cauchy sequence in $X$, consequently we may set
	\[
		Ax := \lim_{j \to \infty} z_j,
	\]
so that we defined a linear operator $A\colon X \to X$ and it satisfies
	\[
		\sigma_j(Ax) = \sigma_j\left( \lim_{k \to \infty} z_k \right)
						   = \lim_{{k \to \infty} \atop {k \geqslant j} } \sigma_j(z_k - z_j)+\sigma_j(z_j)
						   = \sigma_j(z_j)
						   =\left( A_j \circ \sigma_j \right)(x),		
	\]
that is, $A$ turns \eqref{diagram:D2} into a commutative diagram, for every $j$.

Clearly, $Ax$ is well defined, in the sense that it does not depend on the choice of the sequence $(z_j)_j$ in $X$, chosen such that $\sigma_j(z_j) = A_j \circ \sigma_j(x)$  for every $j\in \N$.

If $B\colon X \to X$ is a linear operator such that the diagram \eqref{diagram:D2} is commutative as well, for every $j$, then for every $x \in X$ we have $\sigma_j(Ax-Bx) = [0]_j$ for every $j$, or equivalently, $Ax-Bx \in p_j^{-1}(\{0\})$ for every $j$. Since $(p_j)_j$ is a separating family of seminorms, we get $Ax-Bx=0$, for every $x \in X$. Thus, $A$ is the unique linear operator that turns \eqref{diagram:D2} into a commutative diagram, for every $j$.

We could prove that $A$ is a closed operator and then by the Closed Graph Theorem we would conclude that $A$ is a bounded linear operator on $X$. Instead, we prove a stronger statement.


\medskip
\textbf{\underline{Claim 3:}} $A$ is a strongly compatible operator on $X = \big( X, (p_j)_j \big)$.
\medskip

If $p_j(x)=0$ then $\sigma_j(x)=[0]_j$ and $\sigma_j(Ax) = A_j \circ \sigma_j(x) = [0]_j$, that is, $p_j(Ax)=0$.

Moreover, for every $j$, we see that
	\begin{align*}
		\sup_{p_j(x) \leqslant 1} p_j(Ax) &= \sup_{p_j(x) \leqslant 1}\left \{ \inf_{p_j(z)=0} p_j(Ax) - p_j(z) \right \}\\
													  & \leqslant \sup_{p_j(x) \leqslant 1}\left \{ \inf_{p_j(z)=0} p_j(Ax - z)\right \}
													  = \sup_{p_j(x) \leqslant 1} \| [Ax]_j \|_j \\
													  & \leqslant \sup_{\| [x]_j \|_j \leqslant 1} \| A_j [x]_j \|_j < \infty
	\end{align*}
and the last inequality holds because $\| [x]_j \|_j = \inf_{p_j(y)=0} p_j(x-y) \leqslant p_j(x) \leqslant 1$ whenever $p_j(x) \leqslant 1$.

Therefore, $A \in \Lsc$.

In order to see that $A$ is in fact the infinitesimal generator of $T(\cdot)$, we need to recognize that these projections $\sigma_j$ actually preserve a handy property of projections on Euclidean spaces.

\medskip
\textbf{\underline{Claim 4:}} If $(x_{\lambda})_{\lambda \in \Lambda}$ is a net in $X$ with the property that $[x_{\lambda}]_j \ {\underset{X_j}{\overset{\lambda \in \Lambda}{\longrightarrow}}} \ [0]_j$ for every $j$, then $(x_{\lambda})_{\lambda \in \Lambda}$ is convergent in $X$ and $x_{\lambda} \ {\overset{\lambda \in \Lambda}{\longrightarrow}} \ 0$.
\medskip

We just have to prove that $(x_{\lambda})_{\lambda \in \Lambda}$ is a Cauchy net in $X$. Indeed, given $\varepsilon>0$, we get
	\begin{align*}
		p_j(x_{\lambda} - x_{\eta}) & = \inf_{p_j(z) =0}\left \{ p_j(x_{\lambda} - x_{\eta}) - p_j(z) \right \} \\
												& \leqslant \inf_{p_j(z) =0} p_j(x_{\lambda} - x_{\eta} - z) \\
												& = \| [x_{\lambda}]_j - [x_{\eta}]_j \|_j < \varepsilon,
	\end{align*}
whenever $\lambda, \eta \succeq \gamma$, for some $\gamma = \gamma(\varepsilon, j) \in \Lambda$.

\medskip
\textbf{\underline{Claim 5:}} $A$ is the infinitesimal generator of $\{T(t): t \in \R\} \subset \Lsc$.
\medskip

Given $x \in X$, for every $j \in \N$ we have
	\begin{align*}
		\left[ Ax - \dfrac{T(t)x-x}{t} \right]_j = [Ax]_j -  \dfrac{[T(t)x]_j - [x]_j}{t}
															= \left (A_j[x]_j - \dfrac{T_j(t)[x]_j - [x]_j}{t} \right)
															\quad {\underset{t \to 0}{\overset{X_j}{\longrightarrow}}}  \quad  [0]_j,
	\end{align*}
by the definitions of $A_j$ and $T_j(\cdot)$.

By Claim $4$, we conclude that $\left( \dfrac{T(t)x-x}{t} \right)_{t \in \R}$ converges in $X$ and
	\[
		\lim_{t \to 0} \dfrac{T(t)x-x}{t} = Ax, \mbox{ for every } x \in X.
	\]

Therefore, by Theorem \ref{Th:generation} and by the uniqueness of the infinitesimal generator on Fr\'{e}chet spaces,
\[
	T(t) = \sum_{n=0}^{\infty} \dfrac{t^n}{n!} A^n = e^{tA}, \mbox{ for every } t \in \R.
\]
\endproof

As the reader may note, the spaces $X_j$ were required to be complete so we could use the standard theory on Banach spaces to get the infinitesimal generator $A_j$ of the uniformly continuous group $T_j(\cdot)$ in $X_j$. More precisely, the Neumann series theorem requires completeness and it was implicitly used in the proof. One could remove such requirement from the theorem, since one could redo the proof dealing with the completion of $X_j$, say $\overline{X_j}$: every projection $\pi_j\colon X_{j+1} \to X_j$ can be extended continuously to the completions and we again achieve a commutative diagram with it, so the arguments above still hold. Fortunately, the seminorms of the Fr\'{e}chet spaces $X$ of interest for solving evolution PDEs have a special local property that turns every $X_j$ into a Banach space. Check Proposition \ref{Th:Xj_Banach}.

Just to emphasize, we actually proved a stronger result, that extends the Theorem $1.2$ of Pazy~\cite{Pazy}.

\begin{theorem}\label{Th:equivalence-generation}
A linear operator $A \colon D(A) \subset X \to X$ is the infinitesimal generator of a uniformly continuous group if and only if $A$ is a strongly compatible operator on $X$.
\end{theorem}

Mathematics is surprisingly wondrous: although it was not our intention to use algebraic arguments in the proof of Therem \ref{Th:generator}, we did! Implicitly, Claim 2 and Claim 3 together are a universal property of the projective limit of the spaces $X_j$. This suggests that the strongly compatible operators are good morphisms in the category of all complete locally convex topological vector spaces, where its morphisms are the bounded linear operators on theses spaces. Consequently, we are dealing with convenient bounded operators for our analytical aims, that is, generation of uniformly continuous groups on Fr\'{e}chet spaces.

Let us briefly introduce the subject below. Readers who wish to learn more about it may see~\cites{Narici,Treves,Rotman}, and this includes ourselves.

\begin{remark}[\textbf{Universal property of projective limits}]

Let ${\mathscr C}$ be the category of the complete TVS and the continuous linear applications between them. We may consider a general directed set $(\Lambda, \preceq)$, but $(\N, \leqslant)$ is sufficient for our purpose.

In the proof, we have a projective system; that is, we have a map $\N \ni j \to X_j$, where $X_j$ is a Banach space, and continuous linear maps $\pi_{k,j} \colon X_k \to X_j$, for $k \geqslant j$, with the property that $\pi_{k,j}=\mbox{Id}_{X_k}$ if $k=j$; and $\pi_{k,j} \circ \pi_{l,k} = \pi_{l,j}$ for $l \geqslant k \geqslant j$ in $\N$. In other words, we naturally have a functor $\Phi \colon (\N, \leqslant) \to {\mathscr C}$ defined by
	\[
		 j \quad {\overset{\Phi}{\longmapsto}} \quad X_j \ \mbox{(objects)}
	\]
	and
	\[
		k {\overset{\geqslant}{\longrightarrow}} j \quad {\overset{\Phi}{\longmapsto}} \quad X_k {\overset{\pi_{k,j}}{\longrightarrow}} X_j \ \mbox{(morphisms)}.
	\]

The \textbf{projective limit} of $\Phi$ is an object in ${\mathscr C}$, denoted by
	\[
		\mbox{ proj lim}_{j \in \N} \Phi(j) \mbox{ or simply } \varprojlim_{j \in \N} X_j,
	\]
and a family of continuous linear maps (which are morphisms in ${\mathscr C}$) $\Phi_k \colon \ds\varprojlim_{j \in \N} X_j \to X_k$, $k \in \N$, with the following two properties:
\begin{itemize}
	\item[\textbf{P1)}] if $l \geqslant k$ in $\N$ then
		\[
			\begin{split}
				\xymatrix{
				& \ds\varprojlim_{j \in \N} X_j \ar[dr]^{\Phi_l} \ar[dl]_{\Phi_k} & \\
				X_k & & X_{l} \ar[ll]_{\pi_{l,k}}
				}
			\end{split}
		\]
	is a commutative diagram; and
	
	\item[\textbf{P2)}] \textbf{[Universal property]} If $W$ is an object in ${\mathscr C}$ and $\{B_k \colon W \to X_k\}_{k \in \N}$ is a family of continuous linear maps such that $B_l = B_k \circ \pi_{l,k}$ for $l \geqslant k$, then there exists a unique continuous linear map
		\[
			B \colon W \to \ds\varprojlim_{j \in \N} X_j
		\]
	which turns
		
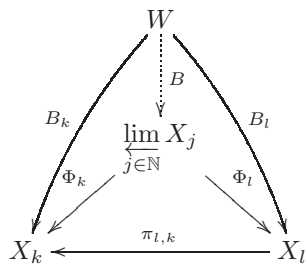
\[
			\begin{split}
				\xymatrix{
				& W \ar@{.>}[d]^{B} \ar@/^/[ddr]^{B_l} \ar@/_/[ddl]_{B_k} & \\
				& \ds\varprojlim_{j \in \N} X_j \ar[dr]^{\Phi_l} \ar[dl]_{\Phi_k} & \\
				X_k & & X_{l} \ar[ll]_{\pi_{l,k}}
				}
			\end{split}
		\]
	into a commutative diagram, whenever $l \geqslant k$ in $\N$.
	\end{itemize}
	
The reader may readily recognize the objects and morphisms used in the proof and hence agree that we proved that the universal property of the projective limit holds there. As pointed out by Babalola~\cite{Babalola}, if $X=\big( X, (p_j)_j \big)$ is a Fr\'{e}chet space with $p_j \leqslant p_{j+1}$ for every $j \in \N$, and every $X_j$ is a Banach space, then
	\[
		\ds\varprojlim_{j \in \N} X_j = X.
	\]
\end{remark}

\begin{corollary}
If $T(\cdot)$ is a uniformly continuous group on $X=\big(X, (p_j)_j \big)$ and every $X_j$ is a Banach space then
	\begin{itemize}
		\item[a)] there exists a unique operator $A$ in $\Lsc$ such that $T(t) = e^{t A}$;
		
		\item[b)] the operator $A$ in part b) is the infinitesimal generator of $T(\cdot)$;

		\item[c)] there exists a sequence $(\omega_j)_{j \in \N}$ of nonnegative numbers such that
		\[
			p_j^{X}\big( T(t) \big) \leqslant \exp(\omega_j \, t), \mbox{ for every } t\in \R;
		\]
		and
		
		\item[d)] the map $\R \ni t \mapsto T(t) \in \Lsc$ is differentiable (in the Fr\'{e}chet sense) and
			\[
				\dfrac{dT(t)}{dt} = A \circ T(t) = T(t) \circ A, \mbox{ for every } t.
			\]
	\end{itemize}
\end{corollary}

The proof is a simple exercise, as long as one has already dealt with such result for Banach spaces.

\begin{remark} 
It is important to observe that the Cauchy problem
	\[
		\left \{
		\begin{array}{l}
		T'(t)=AT(t),\; t\in \R\\
		T(0)=I
		\end{array}
		\right .
	\]
possess a unique solution for each $A\in \Lsc$, by Gronwall's inequality.
\end{remark}

It is noteworthy that if $X$ is a Banach space then the statements of the results and their proofs given above are trivially reduced to the usual proofs. Therefore, the results collected in this section naturally extend the usual theory of generation of uniformly continuous groups on Banach spaces to Fr\'{e}chet spaces. See Pazy~\cite{Pazy}.

At last, we compare the results and aims of this paper with those found in~\cite{Babalola}.

\begin{remark}\label{Remark:comparing-to-Babalola}
We shall resctrict attention to Fr\'{e}chet spaces, although Babalola deals with complete Hausdorff locally convex TVS over $\C$. So let $Y=\big( Y, (q_j){j \in \N} \big)$ be a Fr\'{e}chet space with $q_j \leqslant q_{j+1}$, for $j \in \N$.

A $C_0$-semigroup $\{S(s): s \geqslant 0\}$ in $Y$ (defined as usual) is called a $(C_0,1)$-semigroup if, for every $j \in \N$ and $\delta>0$, there exists a positive constant $c=c(j, S, \delta)$ such that $S(s)V_j \subset c V_j$, for every $0 \leqslant s \leqslant \delta$. Hence, essentially, the author requires strong compatibility over the semigroup with this uniform property on $\delta$. We apply the strong compatibility on the infinitesimal generator instead, and later on, by Theorem~\ref{Th:equivalence-generation}, we recognize that it is equivalent to apply on the (uniformly bounded) group generated by it.

According to~\cite{Babalola}, $(C_0,1)$-semigroups can be characterized as $C_0$-semigroups $\{S(s): s \geqslant 0\}$ on $Y$ such that for each $j \in \N$ there exist a positive number $\sigma_j$ and a natural number $k=k(j)$ so that $q_j\big( S(s)y \big) \leqslant e^{\sigma_j s}q_k(y)$ for every $y \in Y$ and every $s \geqslant 0$. Such result is achieved appealing to arguments of category theory. Actually the ideas used in the proof of Theorem \ref{Th:generator} are closely related to these, although we did not anticipated it. From then on, Babalola's paper is quite different from ours, leading to
	\begin{itemize}
		\item[\textbf{B1)}] the definition of resolvent operators, Hille-Yosida estimates and generation of $(C_0,1)$-semigroups;
		
		\item[\textbf{B2)}] a Trotter-Kato result, that is, a theorem concerning the perturbation of infinitesimal generators; and 
		
		\item[\textbf{B3)}] an application to an ODE, which is an artificial and quite simple one.
	\end{itemize}

On the other hand, we do not really deal with $C_0$-groups or $C_0$-semigroups, but we completely characterize the uniformly bounded groups on Fr\'{e}chet spaces and extend the standard theory to these spaces. As we shall see in the next section, the results above may be applied to the distributions space ${\mathscr F}L^{2}_{loc}$, which provide significant consequences to PDEs; for instance,
	\begin{itemize}
		\item[\textbf{1)}] we explain whether the group preserves or not the support of a function; and
		
		\item[\textbf{2)}] we compare the uniformly bounded group generated by the heat operator (according to Theorem~\ref{Th:generationFL2loc}) with the analytical semigroup generated by it, according to~\cite{Henry}.
	\end{itemize}

\end{remark}

\section{\bf{Some applications to PDEs}}\label{Section:applications}

In Osborne's words~\cite{Osborne}, ``a majority of the topological vector spaces used in analysis are Banach spaces. Also, a majority of the remaining spaces are Fr\'{e}chet spaces. In fact, nearly all the spaces routinely used in analysis are one of four types:
Banach spaces, Fr\'{e}chet spaces, LF-spaces, or the dual spaces of Fr\'{e}chet spaces or LF-spaces''.

That is, the category of Fr\'{e}chet spaces is a large one where the standard Functional Analysis works very well; for instance, Hahn-Banach, Banach-Alaoglu, Banach-Steinhauss, Open Mapping, Closed Graph and Krein-Smulian theorems hold on them. Of course, they hold on more general locally convex spaces. But the whole point of this approach is to consider a phase space (more general than Banach spaces) where the notion of exponential of a bounded linear operator makes sense and may be used to obtain the solution of the Cauchy problem \eqref{eq:Cauchy}.

Originally, the main motivation was the distributional aspects involved in a such problem that were not recognized by the usual approach on normed spaces. Doubtlessly, geometric intuition has been a good guide for solving many differential problems by seeking solution on Hilbert and Banach spaces, thanks to results as Lax-Milgram and regularization theorems and Hille-Yosida theorem. Nevertheless we are convinced that it has restrained our understanding of many phenomena, such as the meaning of the solution of the heat equation for negative times, as we shall see.

It turns out that ${\mathscr F}L^{2}_{loc}$ is a good Fr\'{e}chet space to start with, since it consists of special tempered distributions and it herds the properties of the Fourier transform on $L^{2}$. The next theorem points out that constant coefficients $m$-$\Psi$DOs defined on ${\mathscr F}L^{2}_{loc}$ naturally satisfy the strongly compatibility required by the generation Theorem~\ref{Th:generation}, and hence it is worth the effort so far.

	\begin{theorem}\label{Th:generationFL2loc}
		Every constant coefficients $m$-$\Psi$DO $a(D)$ induces a continuous linear map from $\left({\mathscr F}L^{2}_{loc}, (p_{j}^{*})_{j\in \N} \right)$ to itself by setting
		\[
			a(D)[u] := \big[a(D)u\big], \textrm{ for } [u] \in {\mathscr F}L^{2}_{loc},
		\]
		which is strongly compatible with the seminorms $p_{j}^{*}$ on ${\mathscr F}L^{2}_{loc}$ and consequently generates a group on it.
	\end{theorem}

\proof Recall that $p_{j}^{*}\big([u]\big) = \left (\displaystyle\int_{|\xi|\leqslant j} \big| \widehat{[u]}(\xi) \big|^{2} \, d\xi \right)^{1/2}$, for $[u] \in  {\mathscr F}L^{2}_{loc}$ and $j\in \N$.

Clearly $a(D)$ is linear. To see that it is continuous, first note that if $|\xi| \leqslant j$ then
\[
	[a(D)u] \ \widehat{} \ (\xi) = \lim_{ \substack{l \to \infty \\ L^{2}(B(0,j))}} \big( a(D)u_{l} \big) \ \widehat{} \ (\xi)
											= \lim_{ \substack{l \to \infty \\ L^{2}( B(0,j))}} a(\xi) \, \widehat{u_{l}}(\xi)
											= a(\xi) \, \widehat{[u]}(\xi),
\]
where $\widehat{[u]}$ is the limit in $L^{2}\big( B(0,j) \big)$ of the Fourier transform of some sequence $(u_l)_{l \in \N} \in [u]$, by definition. Now, given $[u] \in {\mathscr F}L^{2}_{loc}$, we have
\begin{align*}
	p_{j}^{*}\big(a(D)[u]\big) &= p_{j}^{*}\Big( \big[a(D)u \big] \Big)
											  = \left (\int_{|\xi|\leqslant j} |a(\xi)|^{2} \, \big| \widehat{[u]}(\xi) \big|^{2} \, d\xi \right)^{1/2} \\
											 & \leqslant \| a \|_{L^{\infty}(B(0,j))} \, \big\| \widehat{[u]} \big\|_{L^{2}(B(0,j))}\\
											 & = \| a \|_{L^{\infty}(B(0,j))} \, p_{j}^{*}\big([u]\big),
\end{align*}
whence we also obtain that $a(D)$ is strongly compatible with the seminorms $p_{j}^{*}$. We therefore conclude the proof by Theorem \ref{Th:generation}.

\endproof

With the definitions and results we have presented so far, we can already obtain the main result described in the Abstract.

\begin{corollary}
If $a(D)$ is a constant coefficients $m$-$\Psi$DO in ${\mathscr F}L^{2}_{loc}$ then the Cauchy problem associated to it, namely
		\begin{equation*}
			\left\{\begin{array}{l}
			u_{t} = a(D)u, t \in \R \\
			u(0) = u_0
			\end{array}
			\right. ,
		\end{equation*}
has a unique solution in ${\mathscr F}L^{2}_{loc}$.
\end{corollary}

As we have seen in Remark \ref{Obs:FL2loc}, the space ${\mathscr E}'(\R^{n})$ is a subspace of ${\mathscr F}L^{2}_{loc}$, so we may wonder whether a semigroup $\{e^{ta(D)}:t \geqslant 0\}$ in ${\mathscr F}L^{2}_{loc}$ lets ${\mathscr E}'(\R^{n})$ invariant. If $n=1$ and the symbol $a=a(\xi)$ is a polynomial then the question is solved below. If $z$ is a complex number, we write $z = \Re\,z+ i  \Im\,z$, with $\Re\,z ,\Im\,z \in \R$.

	\begin{theorem}\label{Th:restriction-E-R}
	
		Let $a(D)=\sum_{\alpha=0}^{m}a_{\alpha}D^{\alpha}$ be a linear differential operator with constant coefficients of order $m$, $a(\xi)=\sum_{\alpha=0}^{m}a_{\alpha}\xi^{\alpha}$ its symbol and $\{e^{ta(D)}: t \in \R\}$ the group generated by it on ${\mathscr F}L^{2}_{loc}$, according to Theorem \ref{Th:generationFL2loc}. Then
			\begin{equation}
				e^{ta(D)}\big( {\mathscr E}'(\R) \big) \subset {\mathscr E}'(\R) \textrm{ for all } t\geq 0
			\end{equation}
		if and only if, either $m=1$ and $\Re\;a_{m}=0$ or $m=4k$, for some $k=0,1,2,\dots $ and $\Re\;a_{m}<0$.	
	\end{theorem}

\proof

Recall that ${\mathscr E}' (\R) \subset {\mathscr F}L^{2}_{loc}$. Given $u\in {\mathscr E}'(\R)$ and $t\in \R$, by Paley-Wiener-Schwartz theorem (\cites{HormanderI,Stein,Treves}), $e^{ta(D)}u\in {\mathscr E}'(\R)$ if and only if, $\widehat{e^{ta(D)}u} \colon\R \to \C$ has an analytic extension $V_{(t,u)} \colon \C \to \C$ and there exist constants $C=C_{(t,u)},R=R_{(t,u)}>0$ and $N=N_{(t,u)}\in \N$ such that for all $z\in \C$ we have
	\[
		|V_{(t,u)}(z)|\leq C_{(t,u)}(1+|z|)^{N_{(t,u)}}e^{R_{(t,u)} |\Im\;z|}.
	\]

Note that $\widehat{e^{ta(D)}u}=e^{ta(\xi)}\hat{u}$ for every $u\in {\mathscr F}L^{2}_{loc}$; in particular, for every $u\in {\mathscr E}' (\R)$. And by hypothesis, $\C \ni z \mapsto e^{ta(z)}\in \C$ is an analytic function. If $u\in {\mathscr E}' (\R)$ then $\xi \mapsto \widehat{e^{ta(D)}u}=e^{ta(\xi)}\hat{u}$ has an analytic extension
	\[
		\C \ni z \mapsto V_{(t,u)}(z)= e^{ta(z)}\hat{u}(z).
	\]

Now fix $ t\geq 0$. For $z\in \C$ the following estimates hold
	\[
		e^{t\Re\,a(z)} |\hat{u}(z)|=|\widehat{e^{ta(D)}u}(z)|=|V_{(t,u)}(z)|\leq C_{(t,u)}(1+|z|)^{N_{(t,u)}}e^{R_{(t,u)} |\Im\;z|}
	\]
if and only if, there exist constants $R,c>0$, such that $e^{\Re\,a(z)}\leq ce^{R |\Im\;z|}$ for every $z\in \C$, if and only if there exist $R=R(a),M=M(a)>0$ such that $\Re\,a(z)\leq R |\Im\;z|$, whenever $|z|\geq M$, if and only if there exist $R'=R'(a),M'=M'(a)>$ such that $\Re(a_{m}z^{m})\leq R' |\Im\;z|$, whenever $|z|\geq M'$.

On the other hand, given $z=\xi+i\eta \in \C$, with $\xi,\eta \in \R$, we put  $a_{m}=\alpha+i\beta$, also with $\alpha,\beta \in \R$. For $m\geq 2$ we have (the cases $m=0,1$ are easy as the reader may verify)
	\begin{equation*}
		\Re(a_{m}z^{m})=\left \{ \begin{array}{l}
		\alpha \xi^{m}+\alpha \eta^{m}+q_{m}(\xi,\eta),\; if \;m=4k\; \text{ for some }k\in \N\\
		\alpha\xi^{m}-\beta \eta^{m}+q_{m}(\xi,\eta),\; if \;m=4k+1\; \text{ for some }k\in \N\\
		\alpha \xi^{m}-\alpha \eta^{m}+q_{m}(\xi,\eta),\; if \;m=4k+2\; \text{ for some }k\in \Z_{+}\\
		\alpha\xi^{m}+\beta \eta^{m}+q_{m}(\xi,\eta),\; if \;m=4k+3\; \text{ for some }k\in \Z_{+},
		\end{array}
		\right.
	\end{equation*}
where $q_{m}(\xi,\eta)$ is a real polynomial of degree $m$ but without the powers $\xi^{m}$ and $\eta^{m}$.

To finish the proof we just have to observe that:

{\bf Case $m=4k$:} 
In this case $m$ is even, so $\xi^{m}\geq 0$ and $\eta^{m}\geq 0$ for all $\xi,\eta \in \R$. If $\alpha <0$, for $|(\xi,\eta)|$ large enough, we have
	\[
		\alpha \xi^{m}+\alpha \eta^{m}+q_{m}(\xi,\eta)=(\alpha \xi^{m}+\alpha \eta^{m})\left(1+\frac{q_{m}(\xi,\eta)}{\alpha \xi^{m}+\alpha \eta^{m}}\right )\leq \frac{\alpha}{2}( \xi^{m}+ \eta^{m})<0<|\eta|.
	\]

If $\alpha >0$, given $c>0$, take $|(\xi,\eta)|$ large enough so that
	\[
		\alpha \xi^{m}+\alpha \eta^{m}+q_{m}(\xi,\eta)=(\alpha \xi^{m}+\alpha \eta^{m})\left(1+\frac{q_{m}(\xi,\eta)}{\alpha \xi^{m}+\alpha \eta^{m}}\right )\geq \frac{\alpha}{2}( \xi^{m}+ \eta^{m})>c|\eta|.
	\]

The case $\alpha =0$ is left to the reader.
\\

{\bf Case $m=4k+1$:} Suppose $\beta >0$. For a given $c\in \N$, we fix $\xi \in \R$ and, since $m \geqslant 1$ is odd, we may choose a negative number $\eta<0$ such that
	\[
		\alpha\xi^{m}-\beta \eta^{m}+q_{m}(\xi,\eta)= \eta^{m}\left(\frac{\alpha\xi^{m}}{ \eta^{m}}-\beta +\frac{q_{m}(\xi,\eta)}{ \eta^{m}}  \right) \geq \frac{-\beta}{2}\eta^{m} >c |\eta|.
	\]

Similarly, we prove for other choices of signs of $\alpha$ and $\beta$.
\\

{\bf Case $m=4k+2$:} Since $m$ is even, we get $\xi^{m}\geq 0$ and $\eta^{m}\geq 0$ for all $\xi,\eta \in \R$. If $\alpha <0$ then, given $c\in \N$, for $\eta \in \R$ large enough, we have
	\[
		\alpha \xi^{m}-\alpha \eta^{m}+q_{m}(\xi,\eta)=-\alpha \eta^{m}\left(\frac{\xi^{m}}{-\eta^{m}}+1+ \frac{q_{m}(\xi,\eta)}{-\alpha \eta^{m}}  \right)\geq \frac{-\alpha}{2}\eta^{m} >c |\eta|.
	\]

Analogously for $\alpha >0$, with $\xi ^{m}$ in the place of $\eta^{m}$; and the case $\alpha =0$ is left to the reader.
\\

{\bf Case $m=4k+3$:} Analogously, we have that $m$ is odd and greater than $1$. If $\beta >0$, fix $\xi \in \R$ and, given $c\in \N$, for $\eta>0$ large enough, we have
	\[
		\alpha\xi^{m}+\beta \eta^{m}+q_{m}(\xi,\eta)= \eta^{m}\left(\frac{\alpha\xi^{m}}{ \eta^{m}}+\beta +\frac{q_{m}(\xi,\eta)}{ \eta^{m}}  \right) \geq \frac{\beta}{2}\eta^{m}>c|\eta|.
	\]

For the other signs of $\alpha$ and $\beta$, the proof follows similarly and the theorem is proved.

\endproof

\begin{remark}\label{sg-distribu}
The last theorem provides conditions under which some differential operators generate semigroups on the space of distributions ${\mathscr E}'(\R)$. It is possible to show that if $u\in {\mathscr E}'(\R)$ and $a(\xi)=\sum_{\alpha=0}^{ m}a_{\alpha}\xi^{\alpha}$ is a polynomial which satisfies the conditions of Theorem \ref{Th:restriction-E-R} then the series of pseudodifferential operators
	\[
		\sum_{n=0}^{\infty}\frac{t^{n}a(D)^{n}}{n!}u
	\]
is such that, for every $\phi \in {\mathscr S}(\R)$ with $\widehat{\phi}\in C^{\infty}_{c}(\R)$, the following ``$\bigstar$-weak convergence'' holds
	\[
		\left \langle \sum_{k=0}^{n}\frac{t^{k}a(D)^{k}}{k!}u, \phi \right \rangle \underset{n\to \infty}{\longrightarrow} \left \langle e^{ta(D)}u , \phi \right \rangle.
	\]

This means that the approach we have considered on ${\mathscr F}L^{2}_{loc}$ might be a good attempt to extend theorems about generation of semigroups to more general space of distributions like ${\mathscr S}'(\R^{n})$ or even ${\mathscr D}'(\Omega)$.

\end{remark}


\vspace{0.5cm}

Similarly, since $L^{2}(\R^{n})\subset {\mathscr F}L^{2}_{loc}$, we wonder whether a group $\{e^{ta(D)}: t\in \R\}$ in ${\mathscr F}L^{2}_{loc}$ lets $L^{2}(\R^{n})$ invariant. The next theorem give a complete characterization of the groups which do that.

\begin{theorem}\label{Th:restriction-L2}

	Let $a(D)=\sum_{|\alpha|\leqslant m}a_{\alpha}D^{\alpha}$ be a linear differential operator with constant coefficients on $\R^n$, let $a\colon \R^{n}\to \C$ be its symbol and $\{e^{ta(D)}: t \in \R \}$ the group generated by it on ${\mathscr F}L^{2}_{loc}$ according to Theorem \ref{Th:generationFL2loc}. 
	
	Then
		\begin{equation}\label{condition0}
			e^{ta(D)}\big( L^{2}(\R^{n})\big) \subset L^{2}(\R^{n}) \textrm{ for all } t \geqslant 0
		\end{equation}
	if and only if
		\begin{equation}\label{condition1}
			\sup_{\xi \in \R^{n}} \,e^{t \, \Re\,\, a(\xi)} <\infty \textrm{ for all } t \geqslant 0.
		\end{equation}
		
\end{theorem}

\proof It is easy to see that $\widehat{e^{ta(D)}u}=e^{ta(\xi)}\hat{u}$ for every $u\in {\mathscr F}L^{2}_{loc}$ and that
\[
	\int _{\R^{n}}|e^{ta(\xi)}\hat{u}(\xi)|^{2} \, d\xi = \int _{\R^{n}}e^{2t \, \Re\,a(\xi)}|\hat{u}(\xi)|^{2} \, d\xi, \textrm{ for } u\in L^{2}(\R^{n}),
\]
so that, given $u \in L^{2}(\R^{n})$, $e^{ta(D)}u$ belongs to $L^{2}(\R^{n})$ if and only if $\int e^{2t \, \Re\,a(\xi)}|\hat{u}(\xi)|^{2} \, d\xi$ is finite.

If $(\ref{condition1})$ holds, given $t\geqslant 0$ and $u\in L^{2}(\R^{n})$, let $M=\underset{|\xi|\geqslant R}{\sup} \,e^{2t \, \Re\,a(\xi)} $ then
\[
	\int _{\R^{n}}e^{2t \, \Re\,a(\xi)}|\hat{u}(\xi)|^{2} \, d\xi \leqslant M\int _{\R^{n}}|\hat{u}(\xi)|^{2} \, d\xi <\infty
\]
and $(\ref{condition0})$ is true.

Conversely, let us prove the contrapositive. If $(\ref{condition1})$ does not hold, then there exist $t \geqslant 0$ and a sequence $(\xi_{N})_{N\in \N}$ in $\R^{n}$ such that $|\xi_{N}|\to \infty $ and
\[
	e^{2t \, \Re\,a(\xi_{N})}\geqslant \frac{2^{N}}{N},
\]
for all natural number $N$.

Now we take a countable collection of disjoint balls $B_{N}:=B(\xi_{N};r_{N})$ such that 
\[
	e^{2t \, \Re\,a(\xi)}\geqslant \frac{2^{N}}{2N}, \textrm{ for all } \xi\in B_{N} \textrm{ and } N\in \N.
\]

To do this, fix a real number $r>0$ such that $B(\xi_{N}; r)\cap B(\xi_{M}; r)=\varnothing$, which is possible since $|\xi_{N}| \to \infty$. By the continuity of $\xi \mapsto e^{2t \, \Re\,\, a(\xi)}$ in $\mathbb{R}^{n}$, for each $N \in \N$, let $r_{N}' >0$ be such that $e^{2t \, \Re\,\, a(\xi)}\geqslant 2^{N}/(2N)$ for all $\xi \in B(\xi_{N};r_{N}')$ and hence set $r_N := \min\{r, r_{N}'\}$.

Let $f_{N}$ be defined by
\[
	f_{N}(\xi) := \frac{2^{-N/2}}{\big(m(B_{N})\big)^{1/2}}\chi_{B_{N}}(\xi), \xi \in \mathbb{R}^{n},
\]
then by the Monotone Convergence Theorem the function $f:=\sum f_{N}$ belongs to $L^{2}(\R^{n})$, because
\[
	\int _{\R^{n}} f^{2} (\xi) \, d\xi = \int _{\R^{n}} \left( \sum_{N=1}^{\infty} f_{N}^{2}(\xi) \right) \, d\xi
													 = \sum_{N=1}^{\infty} \int _{\R^{n}} f_{N}^{2}(\xi) \, d\xi
													 = \sum_{N=1}^{\infty} \frac{1}{2^{N}} < \infty.
\]
 
We claim that $e^{ta(D)}f \, \check{}$ does not belong to $L^{2}(\R^{n})$. Indeed, we apply again the Monotone Convergence Theorem to obtain
\[
	\int _{\R^{n}}e^{2t \, \Re\,\, a(\xi)}|f(\xi)|^{2} \, d\xi = \sum_{N=1}^{\infty} \int_{B_{N}}e^{2t \, \Re\,\, a(\xi)} f_{N}^{2}(\xi) \, d\xi
																			  \geqslant \sum_{N=1}^{\infty}\frac{1}{2N} =\infty.
\]

Hence $(\ref{condition0})$ does not hold for $u:=f \, \check{} \in L^{2}(\R^{n})$ and the theorem is proved.

\endproof

\begin{corollary}\label{InvarianceL2}
If $a(D)$ is a constant coefficients $m$-$\Psi$DO and its symbol $\xi \mapsto a(\xi)$ satisfies $\Re\,\, a(\xi) \leqslant 0$ whenever $|\xi|$ is large enough, then $e^{ta(D)}\big( L^{2}(\R^{n})\big) \subset L^{2}(\R^{n})$ for all $t \geqslant 0$.
\end{corollary}

\bigskip
The last theorem says that we may replace Hille-Yosida theorem (\cite{Pazy}) by the condition $(\ref{condition1})$, if we want to show that a linear differential operator with constant coefficients is the infinitesimal generator of a semigroup on $L^{2}(\R^{n})$. One of the great advantages of Theorem \ref{Th:restriction-L2} over Hille-Yosida theorem is that some differential operators may be infinitesimal generator of semigroups (on a Fr\'{e}chet space) without fulfil the spectral conditions of Hille-Yosida theorem (on Banach spaces). 

Now we provide the result that guarantees that normed quotient spaces associated to the Fr\'{e}chet space ${\mathscr F}L^{2}_{loc}$ are complete.

\begin{proposition}\label{Th:Xj_Banach}
If $X = \big( {\mathscr F}L^{2}_{loc}, (p_j^{*})_{j \in \N} \big)$ then every $X_j$ is a Banach space; more precisely, for every $j\in \N$,
	\[
		X_{j}={\mathscr F}L^{2}_{loc}/(p_j^{*})^{-1}(\{0\}) \equiv L^{2}\big( B[0,j] \big),
	\]
where $B[0,j]$ denotes the closed ball of $\R^n$ centred at the origin and with radius $j\in \N$.
\end{proposition}

\proof
By definition, if $[u]\in {\mathscr F}L^{2}_{loc}$ then
	\[
		\big[ [u] \big]_{j} :=\left \{[u]+[v]:p_{j}^{*}([v])=0\right \}=\left \{[f]\in {\mathscr F}L^{2}_{loc}:\widehat{[f]}(\xi)= \widehat{[u]}(\xi)\; a.e. \; |\xi|\leq j \right \}
	\]
and
	\[
		\Big\| \big[ [u] \big]_{j} \Big\|_{j} = \Big\| \widehat{[u]} \Big\|_{L^{2}(B[0,j])}.
	\]

Hence we may identify $\big[ [u] \big]_{j}$ to $\widehat{[u]}\Big|_{B[0,j]}$ and we are done.

\endproof

Finally, we apply such results to some PDEs and recognize some advantages of this theory over the standard one.

\begin{example}[\bf{The heat equation}]\label{Ex.:heat}

Let us consider the operator
	\[
		A:=1-\Delta \colon H^{2}(\R^{n})\subset L^{2}(\R^{n}, \C) \to L^{2}(\R^{n}, \C),
	\]
where $\Delta=\sum_{j=1}^{n}\frac{\partial ^{2}}{\partial x_{j}^{2}}$ is the Laplacian operator in $\R^{n}$.

By Henry~\cite{Henry}, $A$ is a sectorial operator with $\Re\,\sigma(A)>0$, whence $-(1-\Delta)$ generates an analytic semigroup on $L^{2}$ indicated by $\{e^{-At}:t\geq 0\}$. Besides, the fractional power spaces associated to $A$ are the usual Sobolev spaces $H^{s}=H^{s}(\R^n)$, characterized by the Bessel potentials: $H^{s} = \left\{u\in {\mathscr S}'(\R^n): (1+4\pi ^{2}|\xi|^{2})^{s/2}\widehat{u}\in L^{2} \right \}$.

For $t>0$ and $u\in L^{2}$, we have by Sobolev embedding theorem (found in~\cites{Adams,Bergh_Lofstrom,Folland,Henry,Stein}) that
	\[
		e^{-At}u\in \bigcap_{s\in \R}H^{s} \subset C^{\infty}(\R^{n}),
	\]
so that the solutions of the problem 
	\begin{equation}\label{heat-equation}
		\left\{\begin{array}{l}
		u_{t} +\Delta u= u, t >0 \\
		u(0) = u_0\in L^{2}
		\end{array}
		\right. ,
	\end{equation}
in the sense of~\cites{Pazy, Henry}, are $C^{\infty}$ functions in the variable $x\in \R^{n}$, instantaneously for positive time.

On the other hand, the map $\xi  \mapsto a(\xi)=-(1+4\pi ^{2}|\xi|^{2})$ is the symbol of the pseudodifferential operator $a(D)=-(1-\Delta):{\mathscr F}L^{2}_{loc} \to {\mathscr F}L^{2}_{loc}$, whence $L^{2}$ is left invariant (according to Theorem \ref{Th:restriction-L2}) by the group $\{e^{-t(1-\Delta)}:t\in \R\}$ generated by $a(D)$ on ${\mathscr F}L^{2}_{loc}$, according to Theorem \ref{Th:generationFL2loc}.

Finally, if $u\in {\mathscr F}L^{2}_{loc}$ then
	\[
		\widehat{e^{-t(1-\Delta)}u}=e^{ta(\xi)}\hat{u}, \mbox{ for every } t \geqslant 0,
	\]
and as in page $34$ of Henry~\cite{Henry}, if $u\in L^{2}$ then
	\[
		\widehat{e^{-tA}u}=e^{ta(\xi)}\hat{u}, \mbox{ for every } t \geqslant 0,
	\]
so that both semigroups coincide on $L^{2}$; that is, the group generated on ${\mathscr F}L^{2}_{loc}(\R^{n})$ extends the semigroup generated on the Hilbert spaces. Thus the solution of the heat equation on ${\mathscr F}L^{2}_{loc}(\R^{n})$ for all $t \in \R$ extends the standard solution on Hilbert spaces for $t \geqslant 0$. And hence we are able to solve the heat equation \eqref{heat-equation} {\bf backwards in time} for any initial data $u_{0}\in L^{2} \subset {\mathscr F}L^{2}_{loc}$.

Essentially, for $u_{0}\in L^{2}(\R^{N})$, the regularity of $e^{-tA}u$ has three stages indexed by the time parameter:
	\begin{itemize}
		\item for $t<0$, $e^{-t(1-\Delta)}u_{0}\in {\mathscr F}L^{2}_{loc}$, that is, the solution backwards belongs to a space of very low regularity;
		\item if $t=0$, there is nothing to add, $u_0$ belongs to $L^2$; and
		\item for $t>0$, $e^{-t(1-\Delta)}u_{0}\in C^{\infty}$, that is, the solution are very regular forwards.
	\end{itemize}
	
	Hence it has been suggested that the change of regularity on $t\in \R$ is quite radical: from the space of distributions ${\mathscr F}L^{2}_{loc}$ on negative times to $L^2$ on $t=0$ and then to $C^{\infty}$ instantaneously for positive time. The exponential factor $e^{-t(1+4\pi ^{2}|\xi|^{2})}$ in $\widehat{e^{-t(1-\Delta)}u} = e^{-t(1+4\pi ^{2}|\xi|^{2})}\widehat{u_{0}}$ explains how the regularity of the solution of Heat Equation responds to the time parameter, since
		\[
			\int _{\R^{N}}e^{-2t(1+4\pi ^{2}|\xi|^{2})}(1+|\xi|)^{2M}d\xi <\infty,\text{ for }t>0 \text{ and any }M\in \N,
		\]
and
	\[
		\underset{|\xi|\to \infty}{\lim}e^{-2t(1+4\pi ^{2}|\xi|^{2})}(1+|\xi|)^{2M} = \infty,\text{ for }t<0 \text{ and any }M\in \N.
	\]
\end{example}

\begin{example}[\bf{The derivative operator on $\R$}]

Consider the Cauchy problem
	\begin{equation*}
		\left\{\begin{array}{l}
			u_t = u_x, t \in \R \\
			u(0)=u_0 \in C^{\infty}
		\end{array}
		\right. 
	\end{equation*}
in the phase space $C^{\infty}(\R,\C)$, which is a Fr\'{e}chet space, as seen before.

We shall analyse the generation of a uniformly group on a dense subspace of $C^{\infty}$, without using the generation theorem. This is quite particular to this equation.

Let $\XC$ be the set of all functions $\phi \in C^{\infty}$ such that, for every $ m \in \Z_+$ and $j \in \N$, there exists a constant $M=M(\phi,m,j)>0$ such that
\[
	\sup_{n \in \N} p_{(m,j)} \left( M^{-n} \dfrac{d^n}{dx^n}\phi\right)
		= \sup_{n \in \N} \sup_{|x| \leqslant j} \left| M^{-n} \dfrac{d^{n+m}}{dx^{n+m}}\phi(x) \right| < \infty.
\]

The index ``exp'' refers to the fact that the exponential of the derivative operator is well defined there.

\begin{proposition}

	If $\phi \in \XC$ then $\phi$ is a real analytic function and $\dfrac{d}{dx}\phi \in \XC$.
	
	Moreover,
	\begin{enumerate}
		\item[a)] $\XC$ is a dense subspace of $C^{\infty}(\R)$;
		
		\item[b)] the partial sums $S_N := \ds\sum_{n=0}^{N} \dfrac{t^n}{n!}\dfrac{d^n}{dx^n}\phi$ converges in $C^{\infty}(\R)$ to a function in $\XC$, for every $\phi \in \XC$ and $t \in \R$; its limits is denoted by $e^{t \frac{d}{dx}}\phi$;
		
		\item[c)] $e^{t \frac{d}{dx}} \colon \XC \to \XC$ is well defined and is a bounded linear operator and hence $e^{t \frac{d}{dx}} \in {\mathscr L}\big( C^{\infty}(\R) \big)$ by density;
		
		\item[d)] the family of operators $\{e^{t \frac{d}{dx}}: t \in \R\}$ is a linear uniformly continuous group on $C^{\infty}(\R)$ such that
		\[
			\left( e^{t \frac{d}{dx}}\phi \right)(s) = \phi(s+t), \textrm{ for every } s \in \R.
		\]

	\end{enumerate}
\end{proposition}

\proof

By definition, for every $t \in \R$ and $s_0 \in \supp \phi$, we have $\sup_{|s| \leqslant j} \big| \phi^{(n)}(s) \big| \leqslant M(\phi, j)^{n+1}$, for every $n$, and
	\[
		\left( e^{t \, d/dx} \phi \right)(s_0) = \sum_{n=0}^{\infty} \dfrac{t^n}{n!} \phi^{(n)}(s_0),
	\]
so $\left| \dfrac{1}{n!} \phi^{(n)}(s_0) \right|^{1/n} \leqslant \left( \dfrac{M^{n+1}}{n!} \right)^{1/n} \goesto{\R}{n} 0$. Hence, for every $s \in \supp \phi$,
	\begin{equation}\label{eq:Translation}
		\left( e^{t \, d/dx} \phi \right)(s) = \sum_{n=0}^{\infty} \dfrac{t^n}{n!} \phi^{(n)}(s) = \phi(s+t)
	\end{equation}
and the series has an infinite convergence radius.

Also, if $|x| \leqslant j$ then
	\[
		\left| \dfrac{d^{n+m}}{dx^{n+m}}\phi'(x) \right| = M(\phi, m+1, j)^n  \left| M(\phi, m+1, j)^{-n} \phi^{n+m+1}(x) \right|
				\leqslant M^n \, \sup_{n} p_{(m+1,j)} \big( M^{-n} \phi^{(n)} \big),
	\]
whence $\phi'$ belongs to $\XC$.

\medskip
To prove item a), first note that $x \mapsto e^{-x^2}$ belongs to $\XC$, by choosing $M(m,j) = 2j$. Now we claim that if $\phi \in \XC$ and $f \in C_{c}^{\infty}(\R)$ then $\phi \ast f \in \XC$. Indeed, for $|x| \leqslant j$, we get
	\[
			\left| \dfrac{d^{n+m}}{dx^{n+m}}\big(\phi \ast f \big)(x) \right| \leqslant \int_{B(0,R)} \left| \dfrac{d^{n+m}}{dy^{n+m}} \phi(y) \right| |f(x-y)| \, dy = c(\phi, f, m, j) \, M(\phi, f, m, j)^n,
	\]
where the integer $R=R(f,j)$ is chosen so that $R \geqslant d(0, \supp f) +j$.

Recall that if $\psi \in L^1(\R,\C)$ with $\int \psi = a$, and $f \in L^{\infty}(\R,\C)$ is continuous on an open set $U \subset \R$, then $f \ast \psi_t {\underset{\R}{\overset{t \to 0}{\longrightarrow}}} a \, f$ uniformly on compact sets of $U$, where $\psi_t(x) = t^{-1} \psi(t^{-1}x)$, for $t >0$. See Follland~\cite{Folland}.

In particular, set $\psi(x) := \pi^{-1/2} e^{-x^2}$ and $f \in C_{c}^{\infty}(\R)$. For every pair $(m,j)$ we have
	\[
		p_{(m,j)}(f \ast \psi_t - f) = \sup_{|x| \leqslant j} \big| f^{(m)} \ast \psi_t(x) - f^{(m)}(x) \big| \
												{\underset{\R}{\overset{t \to 0}{\longrightarrow}}} \ 0,
	\]
and then every $f \in C_{c}^{\infty}(\R)$ can be approximated by functions  $f \ast \psi_t \in \XC$ in the topology of $\big(C^{\infty}(\R), (p_{(m,j)}) \big)$. Since $C_{c}^{\infty}(\R)$ is a dense subspace of $C^{\infty}(\R)$, we are done.

\medskip

To prove item b), we just have to verify that $\big( p_{(m,j)}(S_N) \big)_{N \in \N}$ is a Cauchy sequence in $\R$, for every pair $(m,j)$. Indeed, for $M>N$,
	\[
		p_{(m,j)}(S_M-S_N) \leqslant \sum_{n=N+1}^{M} \dfrac{t^n}{n!} p_{(m,j)}\big(\phi^{(m)}\big)
										\leqslant \sum_{n=N+1}^{M} \dfrac{(tM)^n}{n!} \, \sup_n M^{-n} \, p_{(m,j)}\big(\phi^{(m)}\big)
										\ {\underset{\R}{\overset{N,M \to \infty}{\longrightarrow}}} \ 0.
	\]

\medskip

Finally, if $\phi \in \XC$ and $|x| \leqslant j$, we have
	\[
		\left| \dfrac{d^{n+m}}{dx^{n+m}}\Big(e^{t \, d/dx}\phi \Big)(x) \right| \leqslant \sum_{k=0}^{\infty} \dfrac{t^k}{k!} c(\phi, m,j) M^{k+n} = c \, M^n \, e^{t \, M},
	\]
so that $e^{t \, d/dx}\phi \in \XC$.

Clearly, $\phi \mapsto e^{t \, d/dx}\phi$ is linear and, for $(m,j) \in \Z_{+} \times \N$, we have
	\[
		p_{(m,j)} \Big(e^{t \, d/dx}\phi\Big) = \sup_{|x|\leqslant j} \left| \dfrac{d^m}{dx^m} \phi(t+x) \right|
															\leqslant \sup_{|x|\leqslant j + \lceil |t| \rceil} \left| \dfrac{d^m}{dx^m} \phi(t+x) \right|
															= p_{(m, j + \lceil |t| \rceil)}(\phi),
	\]
where $\lceil |t| \rceil$ denotes the smallest integer greater than $|t|$.

Therefore, by \eqref{eq:Translation}, the last two items are proved at once.

\endproof

\end{example}

\begin{example}[\bf{The $i$ derivative operator on $\R$}]
The operator
	\[
		i\frac{d}{dx}: H^{1} \subset L^{2}(\R,\C) \to L^{2}(\R,\C),
	\]
does not fulfil the spectral conditions of Hille-Yosida Theorem, because $\R \subset \sigma(id/dx)$; consequently it cannot generate a semigroup on $L^{2}$ in the sense of this theorem.

However, $i\frac{d}{dx}$ is a pseudodifferential operator and $a(\xi)=-2\pi \xi$ is its symbol, so by Theorem \ref{Th:restriction-L2} we obtain the semigroup $\{e^{itd/dx}:t\geq 0 \}$ on $L^{2}(\R)$.
\end{example}

\begin{example}[\textbf{The bi-Laplacian operator on $\R$}]
By Theorem \ref{Th:restriction-E-R}, the operator $-\frac{d^{4}}{dx^{4}}$ has the property that $e^{-t\frac{d^{4}}{dx^{4}}}u\in {\mathscr E}' (\R)$ whenever $u\in  {\mathscr E}' (\R)$ and $t\geq 0$.
\end{example}

\section{Final Comments}\label{Section:final}

It is known that if $X$ is a Banach space and $A\colon X \to X$ is a bounded linear operator then the exponential of $A$ is well defined as a bounded linear operator, $\exp(A)\colon X \to X$, and may be used to solve the Cauchy problem
	\[
		\left\{\begin{array}{l}
		u_{t} = Au, t \in \R \\
		u(0) = u_0
		\end{array}
		\right. .
	\]

Essentially, if $X$ is a Fr\'{e}chet space, we have recognized which strong compatibility a bounded operator $A$ is required to have so that such resolution still works. Moreover, the approach presented extends naturally the standard theory of generation of uniformly continuous groups and may be applied to pseudodifferential operators with constant coefficients $A=a(D)$ defined in a Fr\'{e}chet space of distributions, namely ${\mathscr F}L^{2}_{loc}(\R^{n})$, which contains $L^2$ and ${\mathscr E}'$. We have established criteria to identify whether the semigroup generated by $a(D)$ acts on these subspaces and analysed the regularization of initial data backwards and forwards by the solution group of the heat equation on ${\mathscr F}L^{2}_{loc}$, which coincides with the standard solution semigroup on Hilbert spaces for positive times.

The strong connection with the usual theory on normed spaces and the results achieved have convinced us that we may consider hyperbolicity (see \cite{aragao1}), semilinear problems, generation of analytic semigroups, non-autonomous linear operators (that is, to consider linear operators $A=A(t)$ depending on $t$) and spectral theory for Fr\'{e}chet spaces as well.

Moreover, Remark \ref{sg-distribu} suggests our approach leads to a new theory of generation of (semi)groups on spaces of distributions and then to the solution of a larger class of linear evolution problems; also, for those problems we already have solved in Banach spaces, it allows us to set a much more singular initial data.

As seen in Example \ref{Ex.:heat}, the theory we presented already explains partially the regularization process which the exponential of the Laplacian operator performs; and it is partial in the sense that it only recognizes three stages of this process. To be more precise, if $u_{0}\in L^{2}$ then $e^{-t(1-\Delta)}u_{0}\in {\mathscr F}L^{2}_{loc}$ for $t<0$, that is, the backwards solution belongs to a space of very low regularity, and this is everything we can say about it for now; $e^{-t(1-\Delta)}u_{0}=u_{0}\in L^{2}$ for $t=0$; and $e^{-t(1-\Delta)}u_{0}\in C^{\infty}$ for $t>0$. We believe that a complete explanation for this phenomena might be provided by a theory of generation of semigroups on spaces of distributions. We mean, if we identify, for every $t\in \R$, which space of distribution the object $e^{-t(1-\Delta)}u_{0}$ belongs to, then the phenomena mentioned before will be satisfactorily explained.

Besides, we know that the inclusion of $\big( L^2, \| \cdot \|_{L^2} \big)$ into ${\mathscr F}L^{2}_{loc}$ is continuous thanks to the Plancherel theorem. As for ${\mathscr E}'$, it is not clear how its original topology is related to its topology as subspace of ${\mathscr F}L^{2}_{loc}$.

Our future aims concerns all these subjects.

\section*{Acknowledgements}

We are deeply thankful for the guidance, support and partnership provided by professor Alexandre Nolasco de Carvalho (University of S\~{a}o Paulo) and professor Walther Hans-Otto (Giessen University), who have inspired us and led us to wonder about consequences and applications of the results we sought.

\newpage
\bibliography{Frechet-10}

\end{document}